\documentclass{amsart}
\usepackage{amssymb,curves}
\newtheorem{Lemma}{Lemma}[section]\newcommand{\bel}{\begin{Lemma}}\newcommand{\eel}{\end{Lemma}}
\newtheorem{Proposition}[Lemma]{Proposition}\newcommand{\bprop}{\begin{Proposition}}\newcommand{\eprop}{\end{Proposition}}
\newtheorem{Theorem}[Lemma]{Theorem}\newcommand{\bthe}{\begin{Theorem}}\newcommand{\ethe}{\end{Theorem}}
\newcommand{\bpr}{{\it Proof.~}}\newcommand{\epr}{$\blacksquare$\\}
\newtheorem{Remark}[Lemma]{Remark}\newcommand{\beR}{\begin{Remark}\rm}\newcommand{\eeR}{\end{Remark}}
\newtheorem{Definition}[Lemma]{Definition}\newcommand{\bed}{\begin{Definition}}\newcommand{\eed}{\end{Definition}}
\newtheorem{Example}[Lemma]{Example}\newcommand{\bex}{\begin{Example}\rm}\newcommand{\eex}{\end{Example}}
\newtheorem{Corollary}[Lemma]{Corollary}\newcommand{\bcor}{\begin{Corollary}\rm}\newcommand{\ecor}{\end{Corollary}}
\newtheorem{Fact}[Lemma]{Fact}\newcommand{\bfact}{\begin{Fact}\rm}\newcommand{\efact}{\end{Fact}}

\newcommand{\beq}{\begin{equation}}\newcommand{\eeq}{\end{equation}}
\newcommand{\bem}{\begin{displaymath}}\newcommand{\eem}{\end{displaymath}}
\newcommand{\beqa}{\begin{eqnarray}}\newcommand{\eeqa}{\end{eqnarray}}
\newcommand{\bee}{\begin{enumerate}}\newcommand{\eee}{\end{enumerate}}
\newcommand{\bei}{\begin{itemize}}\newcommand{\eei}{\end{itemize}}
\newcommand{\bet}{\begin{tabular}{cccccccc}}\newcommand{\eet}{\end{tabular}}
\newcommand{\bpm}{\begin{pmatrix}}\newcommand{\epm}{\end{pmatrix}}
\newcommand{\bM}{\begin{matrix}}\newcommand{\eM}{\end{matrix}}
\newcommand{\ber}{\begin{array}{l}}\newcommand{\eer}{\end{array}}

\newcommand{\tinyM}{\scriptstyle}
\newcommand{\tinyT}{\scriptsize}
\newcommand{\tinyA}{\tinyM\text{\tinyT}}

\newcommand{\cO}{{\mathcal{O}}}

\newcommand{\mC}{\mathbb{C}}
\newcommand{\mN}{\mathbb{N}}
\newcommand{\mP}{\mathbb{P}}\newcommand{\mS}{\mathbb{S}}
\newcommand{\mR}{\mathbb{R}}

\newcommand{\mZ}{\mathbb{Z}}

\newcommand{\al}{\alpha}\newcommand{\be}{\beta}\newcommand{\bc}{\bar\chi}\newcommand{\De}{\Delta}\newcommand{\de}{\delta}
\newcommand{\Ga}{\Gamma}
\newcommand{\ep}{\epsilon}\newcommand{\la}{\lambda}\newcommand{\si}{\sigma}
\newcommand{\Si}{\Sigma}
\newcommand{\ze}{\zeta}

\newcommand{\tE}{\tilde{E}}

\renewcommand{\tt}{\tilde{t}}\newcommand{\tT}{\tilde{T}}\newcommand{\tV}{\tilde{V}}

\newcommand{\lSi}{{\overline{\Sigma}}}

\newcommand{\tBox}{\tilde{\Box}}

\newcommand{\li}{~\\ $\bullet$ }
\newcommand{\ls}{~\\ $\star$ }

\newcommand{\di}{\partial}

\newcommand{\smin}{{\setminus}}\newcommand{\sset}{{\subset}}

\newcommand{\bin}[2]{{#1\choose{#2}}}  

\newcommand{\Nd}{Newton-degenerate }
\newcommand{\dNnd}{directionally Newton non-degenerate}\newcommand{\DNnd}{Directionally Newton non-degenerate}
\newcommand{\tNnd}{topologically Newton non-degenerate} \newcommand{\TNnd}{Topologically Newton non-degenerate}
\newcommand{\gNnd}{generalized Newton non-degenerate}\newcommand{\Nnd}{Newton non-degenerate}

\newcommand{\ND}{Newton diagram}

  \newcommand{\omp}{ordinary multiple point}

\renewcommand{\empty}{\varnothing}
\newcommand{\suml}{\sum\limits}\newcommand{\cupl}{\mathop\cup\limits}\newcommand{\capl}{\mathop\cap\limits}
\newcommand{\prodl}{\prod\limits}\newcommand{\bigcapl}{\mathop\bigcap\limits}

\title{O\MakeLowercase{n some generalizations of} N\MakeLowercase{ewton non-degeneracy for hypersurface singularities}}
\author[]{D\MakeLowercase{mitry} K\MakeLowercase{erner}\\\small B\MakeLowercase{en} G\MakeLowercase{urion}
U\MakeLowercase{niversity}, I\MakeLowercase{srael}}\date{\today}

\address{Department of Mathematics, Ben Gurion University of the Negev, P.O.B. 653, Be'er Sheva 84105, Israel.}
\email{kernerdm@math.bgu.ac.il}

\thanks{The research was constantly supported by the Skirball postdoctoral
fellowship of the Center of Advanced Studies in Mathematics
(Mathematics Department of Ben Gurion University, Israel).
\\
Part of the work was done in Mathematische Forschungsinstitute
Oberwolfach, during the author's stay as an OWL-fellow. Some of the
results were published in the preprint \cite{Kerner-OWP}\date{\today}}


\begin{document}\setcounter{secnumdepth}{6} \setcounter{tocdepth}{1}
\maketitle
\begin{abstract}
We introduce two generalizations of \Nnd~ singularities of
hypersurfaces. Roughly speaking, an isolated hypersurface
singularity is called {\it \tNnd~} if the local embedded topological
singularity type can be restored from a {\it collection} of \ND s (for some coordinate choices).
A singularity that is not \tNnd~is called essentially Newton-degenerate. For plane
curves we give an explicit  characterization of
\tNnd~singularities; for hypersurfaces we provide several examples.

Next, we treat the question: whether \Nnd~or \tNnd~is a property of singularity types or of particular representatives.
Namely, is the non-degeneracy preserved in an equi-singular family? This fact is proved for curves. For hypersurfaces
we give an example of a \Nnd~hypersurface whose equi-singular deformation consists of essentially \Nd~hypersurfaces.

Finally, we define the {\it \dNnd~germs}, a subclass of  \tNnd ~  ones.
For such singularities the classical formulas for the Milnor number and the zeta function
of the \Nnd~hypersurface are generalized.
\end{abstract}
\tableofcontents
\section{Introduction}
We work with germs of complex algebraic (or locally analytic)
hypersurfaces in $\mC^n$, mostly with isolated singularities. By the
singularity type we always mean the local embedded topological type
of a hypersurface germ. For the standard notions of singularity
theory see \cite{AGLV} and \cite{GLSbook}.
\subsection{}To every germ of singular hypersurface (with fixed local analytic coordinates) the \ND~  $\Ga_f$ is associated.
A germ $V_f=\{f=0\}\subset(\mC^n,0)$ is called
\Nnd~(or non-degenerate with respect to its Newton diagram $\Ga_f$) if for each face $\si\in\Ga_f$,
the truncation $f_\si$ of $f$ to $\si$ is non-degenerate
(i.e. the corresponding hypersurface has no singular points in the torus $(\mC^*)^n$). A germ is called {\textit{\gNnd}~if
it is \Nnd~for some choice of coordinates.

The \ND~of a \Nnd~germ is a {\it complete} invariant of the singularity type of the germ.
Namely, if $(V_f,0)$ and $(V_g,0)$ are two \Nnd~germs, such that $\Ga_f=\Ga_g$ then they
have the same embedded topological type \cite{Oka79}.
This distinguishes the \gNnd~germs as especially simple to deal with. For them many
topological invariants of the singularity can be expressed explicitly (or at least estimated)
via the geometry of the \ND~in a relatively simple manner.
For example:
\li the Milnor number \cite{Kouchnirenko76} (cf. also \cite[I.2.1,~pg.122]{GLSbook})
\li the modality (with respect to right equivalence) for functions of two variables (conjectured in \cite[9.9]{Arnol'd74},
proved in \cite[Proposition 7.2]{Kouchnirenko76})
\li the zeta function of monodromy \cite{Varchenko76} (cf. also \cite[II.3.12]{AGLV})
\li the spectrum \cite{Steenbrink76,Khovanski-Varchenko85}(cf. also \cite[II.8.5]{Kulikov98})
\li the Hodge  numbers $h^{p,q}$ \cite{Danilov-Khovanski87}
\li the Bernstein-Saito polynomial \cite{BGMM89}
\li a bound on the \L ojasiewitc invariant \cite{Abderrahmane05} and the \L ojasiewitc-type
inequalities on the sufficiency \cite{Fukui91}.
\\\\
Unfortunately, the condition to be \gNnd~is very restrictive, even in the case of plane curves.
\bex\label{Ex Intro Tangential Decomposition} For the germ $(C,0)\subset(\mC^2,0)$ consider the {\it tangential decomposition}:
$ C=\cupl^k_{i=1}C_\al $.
\\
\parbox{14cm}
{Here each $C_\al $ has unique tangent line $l_i$ (but may contain several branches), and the lines $l_i$ are all distinct.
Therefore, the equation of the tangent cone $T_C$ is $\{l_1^{p_1}\cdots l_k^{p_k}=0\}$, where $p_i$ is
the multiplicity of $(C_\al ,0)$, hence $\sum p_i=p$ \ is
the multiplicity of $(C,0)$. For example, for an \omp\ $p_1=\ldots=p_k=1$.
Note that $p_i=1$ iff $C_\al $ is a smooth branch, not tangent to any other.
The point is that if $(C,0)$ is a \gNnd~
germ then $p_i>1$ for {\it at most} two $i$'s.
}
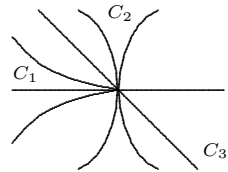
\begin{picture}(0,0)(-60,0)
\curve(-40,0,40,0)\curve(-30,30,30,-30)\put(-40,4){$\tinyA C_1$}\put(32,-25){$\tinyA C_3$}
\put(-4,27){$\tinyA C_2$}
\curve(-40,20,-30,10,-20,5,-10,2,0,0)\curve(-40,-20,-30,-10,-20,-5,-10,-2,0,0)
\curve(-15,30,-5,20,-0.2,0,-5,-20,-15,-30)\curve(15,30,5,20,0.2,0,5,-20,15,-30)
\end{picture}
\\\\
Indeed, first notice that (in some local coordinates) if $(C,0)$ is non-degenerate, then each
$(C_\al ,0)$ is non-degenerate. Moreover, for some fixed  $i$,
if $p_i>1$ and $C_\al $ is non-degenerate with respect to its diagram, then a coordinate axis must
be tangent to $C_\al $ (to reflect the fact that some monomials are absent). Hence, in general, there
are ``not enough coordinate axes" to encode the singularity. In particular,
many singularities with small Milnor numbers and quite simple
defining equations are not \gNnd. For example, the union of three cuspidal branches ($A_2$) with
pairwise distinct tangents are so.
\eex
On the other hand, a locally irreducible plane curve singularity (branch) is \gNnd~ if and only
if it has only one Puiseux pair. An example of not \gNnd~ branch is $(x^2+y^3)^2+x^2y^{3+m}$, $m>0$
with $\mu=15+2m$ (named as $W^\sharp_{1,2m}$ in \cite[\S I.2.3]{AGLV}).

\subsection{Motivation}
The present work has originated from the observation that many germs of curves are ``almost" \gNnd.
Namely, for many of them the singularity type (and thus many of their properties) is reflected on  \ND,
one just has to take {\it several choices} of coordinates.
\bex\label{Ex Intro Collection Newton Diagrams}
Continue the previous example, and
consider again tangential decomposition $C=\cupl_{\alpha=1}^k C_\al $. For each $1\le \alpha\le k$, let $C^{(\alpha)}$
be a germ of curve with the tangential decomposition: $C^{(\alpha)}=\Big(\cupl_{j=1}^{p-p_\alpha}L_j\Big)\cup C_\al $.
Here $\{L_j\}$ are some lines, such that any two are distinct and none is tangent to $C_\al $ (but arbitrary otherwise).
Call such a germ: {\it the directional approximation} of $(C,0)$. (The germ is non-unique, but its topological
singularity type is unique and any two such approximations are connected by a $\mu-$constant family.)
\\
\parbox{14cm}
{If $C_\al $ is \gNnd~then so is each $C^{(\alpha)}$,  and the singularity type of $C_\al $ can be restored from some \ND~of
$C^{(\alpha)}$ (cf. the picture, where the first segment has slope $-1$, the remaining part `is' the diagram of $C_\al $).
Since this diagram is the diagram of the original germ $(C,0)$ too,
the singularity type of $(C,0)$
is completely determined from the {\it collection} of \ND s (corresponding to all the directional approximations).
The precise statement is
in \S\ref{Sec tNnd for Curves}.
}
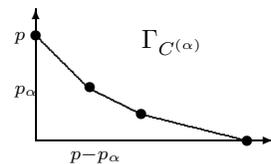
\begin{picture}(0,0)(-20,10)
\put(0,0){\vector(1,0){90}}\put(0,0){\vector(0,1){50}}
\put(18,18){$\bullet$}\put(-2.5,37.5){$\bullet$}\put(37.5,7.5){$\bullet$}\put(77.5,-2.5){$\bullet$}

\curve(0,40,20,20)\curve(20,20,40,10)\curve(40,10,80,0)
\put(-8,38){$\tinyA p$}\put(-8,18){$\tinyA p_\alpha$}\put(13,-7){$\tinyA p-p_\alpha$}  \put(40,35){$ \Ga_{C^{(\alpha)}}$}
\end{picture}
\eex~
Note that if at least one branch of the curve is not \gNnd~then no choice of coordinates can help to recognize
the singularity type in such a way.
\subsection{Results}\label{Sec Intro Results and Questions}
In this work we address some natural questions arising from these examples.
\li Let $(V_f,0)\subset(\mC^n,0)$ be an isolated hypersurface-germ. Suppose all of
its \ND s are known, for any choice of coordinates. {\it Which properties of the germ are determined
by this  information?}

It turns out, e.g. that the collection of all possible diagrams (each labeled with the local coordinate system which determines it)
fix the projectivized tangent cone completely.
Namely, if $f=f_p+f_{p+1}+\ldots$ is the Taylor expansion (in some fixed coordinate system),
 then $f_p$ is determined up
to scaling  from the collection of all diagrams (Proposition \ref{Thm Tangent Cone of tNnd is fixed}).
In fact, even more information can be retrieved (cf. \S \ref{Sec What Is Determined By NDs}), e.g.,
 if $f_p$ is reduced and  $n>2$, then $f_{p+1}$ is also fixed
(up to multiplication of $f$ by an invertible germ).

\li {\it What are the germs whose embedded topological type is completely determined by all the possible \ND s?}
We call such germs: {\it \tNnd} singularities (the name was suggested by E. Shustin). The precise definition is in \S\ref{Sec tNnd}.
It is not clear currently, how broad this class is (in particular all the \gNnd~germs are such), or how to
classify such germs.

For plane curves ($n=2$) we give the complete classification in \S\ref{Sec tNnd for Curves}:
a germ is \tNnd~iff each branch of it is \gNnd~ and the union of any two branches is \gNnd.

For hypersurfaces ($n>2$) the situation is much more complicated.  We give some examples of \tNnd~ germs
in \S\ref{Sec tNnd}. For example all the singularities of Yomdin type
($f_p+f_{p+k}$ with $f_p$ reduced and $f_{p+k}$ generic and the ideal $(f_p,f_{p+k})$ is radical) are such.
Germs that are not \tNnd~are called {\it essentially \Nd}\!.

\li {\it Is  \tNnd~(or \gNnd) a  property of the singularity type or of the germ?}
Namely, suppose an embedded topological
type has a \tNnd~(or \Nnd) realization. Is the generic realization of the type \tNnd~(or \gNnd)?
Or, is this notion preserved in a $\mu-$constant deformation?

This is true for quasi-homogeneous singularities (when the \ND~ is constant along the $\mu-$constant stratum)
by \cite{Varchenko82}, cf. also \cite[III.14.3 Theorem 8]{AGV}. For \Nnd~case the \ND~ can change
essentially along the $\mu-$constant stratum \cite{Briancon-Speder75}. In \cite{Altmann91}
some cohomological conditions on the constancy of \ND~ are given.

Our question can be considered as a weakening of the properties above.
The answer is {\it yes} for the case of curves (Corollary \ref{Thm tNnd and gNnd is a topological property})
and {\it no} for higher dimensions. We give examples in \S\ref{Sec Main Part Germs vs Types} of $\mu^*-$constant
 surface families $(S_t,0)\sset(\mC^3,0)\times\mC^1$ with the central fibre $(S_0,0)$ \Nnd,
 while the generic fibre not \gNnd~ (or even not \tNnd). In fact this situation is typical.

An immediate consequence is the comparison of the equisingular strata versus the ND-topological strata.
Recall that for a fixed \Nnd~representative $(f^{-1}(0),0)\subset(\mC^n,0)$ of a given singularity type, the
ND-topological stratum is defined as the collection of all the hypersurfaces which in some coordinates
 have the diagram $\Ga_f$ and are non-degenerate with respect to $\Ga_f$.
For quasi-homogeneous singularities the ND-topological and $\mu-$constant strata coincide \cite{Varchenko82}.
As our examples show, in many cases the ND-topological strata have positive codimension inside
the equisingular strata.

\li
Once a germ is proven to be \tNnd, its singularity type is determined by the associated collection of \ND s.
Therefore, every topological singularity invariant can be expressed (at least theoretically) via the combinatorics
of the diagrams. However, as the class of \tNnd~singularities is so broad, it seems difficult
to do this generally. Rather, we restrict to a subclass of {\it \dNnd~}singularities (introduced in
 \S\ref{Sec DNnd singularities}).
These are germs with \Nnd~directional approximation (the natural generalization of the case of curves,
example \ref{Ex Intro Collection Newton Diagrams}).
The number of \ND s needed to determine the singularity type is bounded in this case.
For example, if the projectivized tangent cone $\mP T_{(V_f,0)}$ of the hypersurface has only isolated singularities,
then the singularity type of $(V_f,0)$ can be determined by $|Sing(\mP T_{(V_f,0)})|$ coordinate choices.

For \dNnd~singularities we generalize some classical formulas.
In particular in \S\ref{Sec Invariants} the formulas for the Milnor number and
 for the zeta function of monodromy are generalized. For curves we generalize also the formula for modality.

\subsection{Conventions and notations}\label{Sec Intro Conventions and Notations}
In general, we work in the space of the locally analytic hypersurface germs in $\mC^n$. Sometimes we
pass to the space of
germs of (high) bounded degrees (to have a finite dimensional space, to use algebraicity and Zariski topology).
As the singularities are isolated this is always possible by finite determinacy.

Note that for $u$ a locally invertible function: $\Ga_f=\Ga_{uf}$. Hence the \ND~is well defined
by the zero set $V_f$.
The \ND s are assumed to be commode or convenient ($\Ga_f$ intersects all the coordinate axes), unless explicitly stated.
Denote by $f_\si$ the restriction of the function $f$ to the face $\si\in\Ga_f$.
Denote by $\Ga_{f}^-$ the set of real points on or below $\Ga_f$.

Throughout the paper we use equisingular deformations and $\mu-$constant deformations.
In most cases the two notions coincide: for $n\neq3$ this follows by \cite{Le-Ramanujam76},
for $n=3$ and deformations $f_t$ linear in $t$ it follows by \cite{Parusin'ski99}.
In fact our $\mu-$constant deformations are often even $\mu^*-$constant.
\subsection{Acknowledgements}
This work would be impossible without numerous important discussions with G.M. Greuel, A. N\'emethi and E. Shustin.
 A. N\'emethi pointed an uncountable number of errors in the first $N\gg0$ versions of the text.
Many thanks are also to V. Goryunov for important advices.

Part of this work was done during my stay at the Mathematische Forschungsinsitute Oberwolfach (Germany).
 Many thanks to the staff for the excellent working atmosphere.

\section{What is determined by the collection of \ND s?}\label{Sec What Is Determined By NDs}
In this section $(V_f,0)$ and $(V_g,0)\subset(\mC^n,0)$ are germs of isolated hypersurface singularities, such
that for {\em any choice} of coordinates $\Ga_f=\Ga_g$.
Let $f=f_p+f_{p+1}+\ldots$ and $g=g_p+g_{p+1}+\ldots $ be the Taylor expansions in some fixed local coordinate system.
For any $k\geq 0$, let $I_k(f):=\langle f_p,f_{p+1},..,f_{p+k}\rangle\subset\cO:=\cO_{(\mC^n,0)}$ be the ideal generated
by the corresponding homogeneous forms. Finally, let $Rad(I)$ denote the radical of the ideal $I$.
\bprop\label{Thm Tangent Cone of tNnd is fixed}\label{Thm Higher Terms of Taylor Expansion of TNnd are fixed}
Under the above assumptions one has:
\begin{enumerate}
\item  For any $k\geq 0$ the radicals of ideals coincide: $Rad(I_k(f))=Rad(I_k(g))$.
In particular, for $k>0$,  $\cO/Rad(I_{k-1}(f))\supset Rad(\langle [f_{p+k}]\rangle) =Rad(\langle [g_{p+k}]\rangle) \subset\cO/Rad(I_{k-1}(g))$.
\item Moreover, regarding the tangent cone, one has
$f_p=g_p$, up to multiplication by a constant.
\item For a fixed $k>0$, suppose that $I_i(f)$ is radical for any $0\le i\le k$.
Then $I_k(f)=I_k(g)$ and
$\suml_{i=0}^kf_{p+i}=\suml_{i=0}^k(a_i+b_i(x)g_{p+i})$, where $a_i\in\mC^*$ and $b_i(x):(\mC^n,0)\to(\mC,0)$ for all $i$.
\end{enumerate}
\eprop
\bpr (1) Consider the projective hypersurfaces $\{\mP V_{f_{p+i}}\subset\mP^{n-1}\}_{0\leq i\leq k}$. Suppose $x\in\capl^k_{i=0}\mP V_{f_{p+i}}$.
By $GL(n)$ action, one can assume that $x=[0,\ldots,0,1]$. Thus, none of the equations $f_{p+i}$ contains the monomial $x^{p+i}_n$.
In particular, the diagram $\Ga_f$
intersects the $x_n$-axis at a point higher than $(p+k)$. By the assumption on $f$ and $g$,  the same fact holds for $\Ga_g$.
Therefore $x\in\capl^k_{i=0}\mP V_{g_{p+i}}$. Hence, (set theoretically)
$\Big(\capl^k_{i=0}\mP V_{f_{p+i}}\Big)_{red}=\Big(\capl^k_{i=0}\mP V_{g_{p+i}}\Big)_{red}$, and
by Nullstellensatz (1) follows.

(2)
Write $f_p=\prodl_{i=1}^l f_{p,i}^{n_i}$ and
$g_p=\prodl_{i=1}^l g_{p,i}^{m_i}$ for the prime decompositions of $f$ and $g$, where
$f_{p,i}=g_{p,i}$ up to scaling, cf. part (1). Let $x\in V_{f_{p,i}}$ be a
generic point, so that $x$ is a smooth point of the reduced cone $V_{\prod f_{p,i}}$. Apply
linear transformation $\phi$ to $(\mC^n,0)$ in order to get $x=[1,0,\ldots,0]$. Then the monomial $x_1^{deg(f_{p,i})}$ does
not appear in $f_{p,i}$, while for any $j\neq i$ the monomial $x_1^{deg(f_{p,j})}$ does appear in $f_{p,j}$.
Thus, the number $p-deg(f_{p,i})n_{p,i}$ can be restored from the \ND~$\Ga_{\phi(f)}$ by checking the monomial
containing the highest power of $x_1$. Hence, by equality of the \ND s, one gets $n_i=m_i$.
This shows that the scheme structure of the projectivized tangent cone is also restored from the collection of the \ND s.

(3) The coincidence of the ideals is proved by induction. First, note that $\langle f_p\rangle=\langle g_p\rangle$.
Suppose that for all $0\leq i\leq k$ the ideals $I_i(f)$ are radical and
$I_{k-1}(f)=I_{k-1}(g)$.
Hence, by (1) and induction,  $I_k(f)=Rad(I_k(f))=Rad(I_k(g))=Rad(\langle f_p,\ldots,f_{p+k-1},g_{p+k}\rangle)$.
Hence $g_{p+k}=\al f_{p+k}+\sum_{i<k} \be_i f_{p+i}$. Using grading one gets that $\alpha$ is a (zero or non-zero)
constant. Analyzing both cases, $I_k(f)=I_k(g)$ follows.
\epr
\beR The possible attempts to strengthen/generalize the above proposition are obstructed:

(1) It is important to take radicals (for the first statement),
as the following example shows. Suppose $f_p=x^p_1$ with $p\ge3$ and $f_{p+1}$
contains the monomial $x^{p+1}_2$.
 Then all the  relevant coordinate transformations have the form
$x_1\mapsto x_1+Q(x_1,\ldots,x_n)$ with $Q$ quadratic, and the other coordinates $(x_1,\ldots, x_n)$ are moved by
a linear transformation of $GL(\mC^{n-1}_{x_2,\ldots,x_n})$ (which preserves the hyperplane $x_1=0$).
Write $f_{p+1}$ as $\tilde{f}_{p+1}(x_2,\ldots,x_n)+x_1\tilde{f}_p$. Then the possible \ND s
fix the homogeneous form $\tilde{f}_{p+1}$
completely, but impose no restrictions on $\tilde{f}_p$.
E.g., $f=x^p_1+x^{p+1}_2$ and $g=x^p_1+x_1x^p_2+x^{p+1}_2$ have equal \ND s
in any coordinates, but of course $I_1(f)\neq I_1(g)$.

(1') Similarly,  the property
$\cO/Rad(I_{k-1}(f))\supset Rad(\langle [f_{p+k}]\rangle)=Rad(\langle [g_{p+k}]\rangle) \subset\cO/Rad(I_{k-1}(g))$
does not hold without the radicals, as it is exemplified by
the plane curves defined by the polynomials  $f=x^2y^2+x^5+y^5$ and  $g=x^2y^2+x^5-y^5$.
Note that the corresponding \ND s coincide in any coordinate system. (By the direct check, a change
of coordinates whose linear part is identity has no influence on the diagrams.)
But, clearly,
$\cO/Rad(\langle f_4\rangle)\supset \langle[f_5]\rangle\neq\langle [g_5]\rangle \subset\cO/Rad(\langle g_4\rangle)$.

(2) It is not possible to consider the filtration determined
by the \ND~ instead of Taylor expansion.
Indeed, even in the case of a quasi-homogeneous filtration, the lowest order parts do not necessarily coincide.
As an example, consider  $f(x,y,z)=(z^3+x^4y^5)^r+x^{8r}+y^{10r}+\tilde{f}$ and
$g(x,y,z)=z^{3r}+x^{4r}y^{5r}+x^{8r}+y^{10r}+\tilde{f}$.
Here $r\ge2$ and $\tilde{f}$ consists of monomials above the hyperplane $\frac{z}{3r}+\frac{x}{8r}+\frac{y}{10r}=1$.
Note that $f,g$ are \Nnd~ and for any coordinate system $\Ga_f=\Ga_g$. This last statement can be verified as follows.
Take a locally analytic transformation $\phi$ of $(\mC^3,0)$.
If its linear part mixes the coordinates (i.e. if it is not diagonal or a permutation) then obviously
$\Ga_{\phi^*(f)}=\Ga_{\phi^*(g)}$. Hence, one can assume that the linear part is the identity.
Analyzing the non-linear parts, all the  relevant cases are
 $z\mapsto z+\varphi$, where $\varphi \in\langle x^2,xy,xy^2,y^2,y^3\rangle$. By a direct verification,
 one gets again $\Ga_{\phi^*(f)}=\Ga_{\phi^*(g)}$. But the lowest order parts of $f,g$ differ significantly.

(3) (Regarding the last statement of the proposition.) If $I_k$ is radical then $I_0$ is such i.e. $f_p$ is reduced.
However, in general,  the fact that $f_p$ is reduced seems to be not enough to prove that the other intermediate ideals are reduced as well (at least not by general theory of ideals). Below is an example, communicated to
me by D. Eisenbud and B. Ulrich to whom I am very grateful, which shows the subtlety of the problem.

Consider the ideal in $\mC[a,b,c,d,e,f,g]$ generated by the 3 minors of the matrix $\bpm a^2&d^3&f^4\\bc&e^3&g^4\epm$.
The 3 equations $f_5,f_6,f_7$ have degrees 5,6,7. The total ideal is reduced (this can be checked using e.g. \cite{Singular}).
But the ideal generated by the first two equations is non-reduced, e.g. it contains $a^2f_7$ but not $af_7$.

Note that this example is general, it does not use
 all the assumptions of the proposition above. In particular the singularity of this example is non-isolated.
 It would be  interesting to give a counterexample which fit exactly to our situation.

\eeR
\bcor\label{Thm f is Yomdin then g is Yomdin} With the assumptions of proposition
\ref{Thm Tangent Cone of tNnd is fixed} one has:
\begin{enumerate}
\item Let $f=f_p+f_{p+k}+\cdots$ be the Taylor expansion
and suppose that the ideal $\langle f_p,f_{p+k}\rangle$ is radical.
 Then $g$ is contact equivalent  to $f_p+f_{p+k}+${\it some higher order terms} (which might be different from those of $f$).
\item Let $f=f_p+f_{p+k}+f_{p+q}+\cdots$ with $q>k$ be the Taylor expansion and
suppose that the ideal $\langle f_p,f_{p+k},f_{p+q}\rangle$
is radical. Then $g$ is contact equivalent  to $f_p+f_{p+k}+\la f_{p+q}+${\it some higher order terms},
for some $\la\in\mC$.
\end{enumerate}
\ecor
\bpr
We prove the second claim, the first is proved similarly. The last part of proposition
\ref{Thm Tangent Cone of tNnd is fixed} gives: $g=af_p+bf_{p+k}+cf_{p+q}+\cdots$
for $a,b,c$ locally analytic and locally invertible. Hence there exists a locally analytic coordinate scaling:
$(x_1,\cdots,x_n)\to \de(x_1,\cdots,x_n)$ such that $\frac{a}{\de^p}=\frac{b}{\de^{p+k}}$. Apply it and divide $g$ by
$\frac{a}{\de^p}$ to get: $g\sim f_p+f_{p+k}+\la f_{p+q}+..$.
\epr
\beR Regarding  the possible converse of proposition \ref{Thm Higher Terms of Taylor Expansion of TNnd are fixed}
we note the following:
\\
(a) Suppose that the components of the Taylor expansion (for some choice of coordinates) satisfy: $\suml_{i=p}^{p+k}f_i=\suml_{i=p}^{p+k}g_i(a_i+\al_i(x))$
with $a_i\in\mC^*$ and $\al_i:(\mC^n,0)\to(\mC,0)$. Then the \ND s (in that fixed coordinate system)
coincide up to the order $p+k$ (i.e. the parts of $\Ga_f,\Ga_g$ lying below the hyperplane $\sum_i x_i=p+k$ coincide).

This follows immediately from the fact that $\Ga_{\sum_i f_i}=\Ga_{\sum_i a_if_i}$ for $a_i\in\mC^*$ and
$\Ga_{\sum_i f_i}=\Ga_{\sum_i f_i(1+\al_i(x))}$.

On the other hand, the equality $\suml_{i=p}^{p+k}f_i=\suml_{i=p}^{p+k}g_i(a_i+\al_i(x))$ in {\it some}
coordinate system does not imply
that the (truncated) \ND s coincide for {\it any} choice of coordinates. An elementary example is: $f=xy+y^4$ and $g=xy-y^4$.
\\
(b) The equality of radical ideals does not imply any relation of \ND s or singularity types.
As an example consider $f=x^p+xy^p+y^q$ and $g=x^p+y^q$ (for $q>p+1$).
\eeR

\section{\TNnd~hypersurfaces}\label{Sec tNnd}
\subsection{Preparations for the definition.}\label{Sec tNnd preparations for definition}
Start from the following observation. Let $(V_f,0)=\{f=0\}\subset(\mC^n,0)$ be a \gNnd~
isolated singularity. Fix some coordinates $(x_1, \ldots,x_n)$.
Let $\phi\circlearrowright(\mC^n,0)$ be a locally analytic coordinate change, such that $\phi^*(f)$
is non-degenerate with respect to its diagram  $\Ga_{\phi^*(f)}$.
In the space of all the locally analytic series at the origin (or, in some of its truncations if it is necessary)
consider the stratum:
\beq
\Si_{(\phi,\Ga_{\phi^*(f)})}:=\overline{\{g\in\mC\{x\}~|~\Ga_{\phi^*(g)}=\Ga_{\phi^*(f)}\}}.
\eeq
Here the closure is taken in the classical topology (for the coefficients of the defining series).
This stratum is irreducible. Then for the
generic point $g\in\Si_{(\phi,\Ga_{\phi^*(f)})}$
 the {\it local embedded topological types} of $(V_f,0)$ and $(V_g,0)$ coincide \cite{Kouchnirenko76}.

This can be rephrased as follows: any small deformation of $f$ inside $\Si_{(\phi,\Ga_{\phi^*(f)})}$ is $\mu-$constant
(cf. \S\ref{Sec Intro Conventions and Notations}).
\\\\
Recall the notion of Newton weight function \cite[I.3.8]{AGLV} associated to every commode \ND.
Namely, $\la_\Ga:\mR^n_+\to\mR_+$ is defined uniquely by the conditions: $\la_\Ga(\al \vec{x})=\al\la_\Ga(\vec{x})$
 (for any $\al\in\mR_+$) and $\la_\Ga(\Ga)=1$.

Given two diagrams we say $\Ga_1\ge\Ga_2$ if $\la_{\Ga_1}(\Ga_2)\le1$ (or $\la_{\Ga_2}(\Ga_1)\ge1$).
\\\\\\
Suppose a collection of pairs $\{(\phi_i,\Ga_i)\}$ is given  (with $\phi_i\circlearrowright(\mC^n,0)$ local
coordinate changes and $\Ga_i$  \ND s).
\bed
The stratum of hypersurfaces germs, associated with the collection  $\{(\phi_i,\Ga_i)\}$ is the closure of
 the set of all the germs giving the
prescribed diagrams in the prescribed coordinates, i.e.
\beq
\Si_ {\{(\phi_i,\Ga_i)\}}:=\{g\in\mC\{x\}~|~\Ga_{\phi^*_i(g)}\ge\Ga_i~\text{ for all }i\}.
\eeq
\eed
\bel
For any collection $\{(\phi_i,\Ga_i)_i\}$, as above, the associated stratum $\Si_ {\{(\phi_i,\Ga_i)\}}$ is
a (non-trivial) linear subspace of the space of all locally analytic functions at the origin.
In particular it is  closed, irreducible and the
notion of the generic point is well defined.
\eel
\bpr
The condition $\Ga_{\phi^*_i(g)}\ge\Ga_i$ means the absence of some monomials in the Taylor expansion of $\phi^*_i(g)$.
This says that some directional derivatives of $\phi^*_i(g)$ vanish:
$\sum a_{j_1..j_n}\di^{j_1}_{y_1}..\di^{j_n}_{y_n}(g\circ\phi_i)=0$. Here $\{y_j=\phi_i(x_j)\}$ are
the new coordinates. And these conditions are linear in $g$
in any coordinate system.
\epr
\subsection{The main definition.}\label{Sec tNnd the main definition}
\bed\label{Def tNnd}
The function $f\in\mC\{x\}$  is called \tNnd~ if there exist a finite number of coordinate choices
(i.e. locally analytic $\phi_i\circlearrowright(\mC^n,0)$) and the \ND s $\{\Ga_i\}$ such that
$\{\Ga_{\phi^*_i(f)}=\Ga_i\}$ and any small deformation of $f$ inside
 the stratum $\Si_{\{(\phi_i,\Ga_i)\}}$ is $\mu-$constant.
\eed
Recall that for a locally invertible $u\in\mC\{x\}$ one has $\Ga_{uf}=\Ga_f$ and $(V_{uf},0)=(V_{f},0)$.
Moreover, $f$ is \tNnd~iff $uf$ is such. Thus the hypersurface germ $(V,0)\subset(\mC^n,0)$ is defined to
be \tNnd~ if one (and hence any) of its locally defining functions is \tNnd.

The definition \ref{Def tNnd} is equivalent to the following: the general point of $\Si_{\{(\phi_i,\Ga_i)\}}$
corresponds to a hypersurface germ, whose singularity type is that of $(V_f,0)$.
General here means: lying in the complement of a proper analytic subset.
Some comments are in order. Note that the $\mu-$constant deformation is equisingular
  (for $n=3$  use the linearity of the space  $\Si_{\{(\phi_i,\Ga_i)\}}$ and see the remark in
\S\ref{Sec Intro Conventions and Notations}). Thus by semi-continuity of $\mu$ we
deduce that the subset of the germs $f\in \Si_{\{(\phi_i,\Ga_i)\}}$ with the given topological type
is open; and for germs in the complement the Milnor number is strictly larger. So $f$ is \tNnd~iff
it belongs to this (Zariski) open set.

\bex\label{Examples of tNnd}\bee
\item Every \gNnd~ germ is \tNnd. In this case, by definition, just one pair $(\phi,\Ga)$ suffices.
\item Let $C=\cupl_{i=1}^k C_\al $ be the tangential decomposition of a plane curve singularity (cf. example
 \ref{Ex Intro Tangential Decomposition}). If each of $C_\al $ is \gNnd~then $C$ is \tNnd. Indeed, make $k$
 choices of coordinates $(x^{(i)}_1,x^{(i)}_2)$ with $\hat{x}^{(i)}_2$ axis generic and $\hat{x}^{(i)}_1$ axis
 chosen such that the germ $C_\al $ is \Nnd.
 Then get the collection of \ND s similar to those of example \ref{Ex Intro Collection Newton Diagrams}.
 Obviously, this collection of diagrams specifies the singularity type uniquely (the diagram $\Ga_i$ specifies
 the type of $C_\al $ and the fact that the tangent line of $C_\al $ is not tangent to any other $C_{j\neq i}$.

\item The curve germ  $(x^2-y^3)(x^2-y^3+x^3)=0$ is the union of two branches, each being \Nnd,
 but the union is not \tNnd.
It is easy to see that the family $C_{t\neq0}=\{(x^2-y^3)(x^2-ty^3+x^3)=0\}$ has constant \ND~
 for any choice of coordinates. But $\mu(C_{t\neq1})<\mu(C_1)$.
\item \label{Thm SuperIsolated Sings are tNnd} (Yomdin series.)
Let $f=f_p+f_{p+k}+\cdots$ with $f_{p+k}$ generic with respect to $f_p$.
Namely, $\mP V_f$ and $\mP V_{p+k}$ intersect transversally in $\mP^{n-1}$, in particular
$Sing(\mP V_f)\cap \mP V_{p+k}=\empty$.
 Assume that the ideal $\langle f_p,f_{p+k}\rangle$ is radical.
Then $(V_f,0)\subset(\mC^n,0)$ is \tNnd.

Indeed, if for some $g\in\mC\{x\}$ the Newton diagrams $\Ga_g,\Ga_f$  coincide in any coordinates,
then by corollary \ref{Thm f is Yomdin then g is Yomdin} $g$ is contact equivalent with a germ
of the form $f_p+f_{p+k}+$higher order terms. But all these germs are
topologically contact equivalent to $f_p+f_{p+k}$, hence $g$ and $f$ are topologically equivalent.
Finally note that the set of all the coordinate systems can be replaced by a finite set
(by the argument as in the proposition below).

We believe that the condition "$\langle f_p,f_{p+k}\rangle$ is radical" is unnecessary here, but do not
have any rigorous proof.
\eee\eex

\bprop(Consistency of the definition.)\bee
\item Let $(V_f,0)\subset(\mC^n,0)$ be a \tNnd~ germ and $\{(\phi_i,\Ga_i)_{i=1}^k\}$
a collection of pairs fulfilling the condition of definition \ref{Def tNnd}
(i.e. specifying the singularity type of $V_f$
uniquely). Then for any additional pair $(\phi_{k+1},\Ga_{k+1})$  the collection $\{(\phi_i,\Ga_i)_{i=1}^{k+1}\}$
also fulfills the condition of the definition.
\item Let $\{(\phi_i,\Ga_i)_{i}\}_{i\in I}$ be an infinite collection such that
the associated stratum  $\Si_{\{(\phi_i,\Ga_i)\}_{i\in I}}$ (defined similarly to the case of finite collection)
is of positive dimension.
Suppose that any small deformation of $f$ inside this stratum is $\mu-$constant.
Then there exists a finite sub-collection $J\subset I$ such that the conditions of the definition \ref{Def tNnd}
are satisfied for $\{(\phi_i,\Ga_i)_{i\in J}\}$ and hence $(V_f,0)$ is \tNnd.
\eee\eprop
\bpr\bee
\item By assumption any small deformation of $(V_f,0)$ inside the stratum $\Si_{\{(\phi_i,\Ga_i)_{i=1}^k\}}$
is $\mu-$constant,
hence this is true for the substratum $\Si_{\{(\phi_i,\Ga_i)_{i=1}^{k+1}\}}$ as well.
\item
By finite determinacy one can pass to finite jets $J_N$ of some high order.
Let $j_N(g)$ be the $N$-jet of $g$. Then $\mu(V_f,0)=\mu(V_{j_N(f)},0)$
and any small deformation of $(V_{j_N(f)},0)$ inside $\Si_{\{(\phi_i,\Ga_i)\}_{i\in I}}\cap J_N$
is $\mu-$constant.
Present the stratum $\Si_{\{(\phi_i,\Ga_i)\}_{i\in I}}\cap J_N$ as successive intersections:
$\bigcapl_{i\in J}\Si_{(\phi_i,\Ga_i)}\cap J_N$, where $J\subset I$ is an increasing filtration of $I$.

At each step we get a linear subspace of the {\it finite dimensional} space of $N$-jets. Each
intersection either decreases the dimension or has no influence. Therefore the process stabilizes after a finite number
of intersections. Hence, there exists a finite subset $J\subset I$ satisfying the conditions of
the definition.
\epr\eee

\subsection{The case of curves}\label{Sec tNnd for Curves}
For curves it is possible to give a very explicit equivalent definition of topological Newton-non-degeneracy.
\bprop\label{Thm  Criterion on tNnd for curves}
 Let $(C,0)=\cup_i (C_\al ,0)$ be the tangential decomposition. Then $(C,0)$ is \tNnd~ iff each $(C_\al ,0)$ is \tNnd.
Moreover $(C_\al ,0)$ is \tNnd~ iff the following two conditions are satisfied:
\li Each branch of $C_i$ is \gNnd, i.e. locally it is of the type $x^p+y^q$ with $(p,q)=1$.
\li In addition, the union of any two singular branches is a \gNnd~singularity. More precisely, there does
not exist a pair of singular
branches in $C_\al$ with local equation (in some coordinates): $(x^p+y^q+\cdots)(x^p+y^q+\cdots)$. Here the dots mean higher
order terms (i.e. monomials lying over the \ND~of $x^p+y^q$).
\eprop
The last condition can be rephrased as follows: for any pair of singular branches $C_i,C_j$ with the Puiseux pairs $(p_i,q_i)$
and $(p_j,q_j)$ either $(p_i,q_i)\neq(p_j,q_j)$ or $(p_i,q_i)=(p_j,q_j)$ but the intersection multiplicity
$mult_0(C_i\cap C_j)\le p_iq_i$.

\bpr
For the first part note that the collection of $\{\phi_i,\Ga_i\}$ that specifies the type of $(C,0)$
does the same for each $(C_\al,0)$ independently. Conversely, let $\cup_i\{\phi_{\al,i},\Ga_{\al,i}\}$ be
a collection that specifies the type of $(C_\al,0)$. Then the total
collection $\cup_{i,\al}\{\phi_{\al,i},\Ga_{\al,i}\}$ specifies the type of $(C,0)$.

Thus, in the sequel we  assume that all the branches of $(C,0)$ have a common tangent.
\\$\Rrightarrow$ Suppose $(C,0)$ contains a branch $(C_i,0)$ which is not \gNnd.
Choose coordinates in which the defining equation of $C_i$ can be written as: $(x^{ar}+\cdots+y^{br})+h=0$.
Here $1<a<b$, $(a,b)=1$, the part $(x^{ar}+\cdots+y^{br})$ is quasi-homogeneous and degenerate and $h$
is of higher order with respect to the weights above. Let $(C'_i,0)$
be a generic germ with the \ND~of $C_i$, in particular $C'_i$ is \Nnd.
Note, that the \ND s of $C_\al ,C'_i$ coincide in any coordinates.

Consider the germ $C'=\bigcup_{j\neq i}C_j\cup C'_i$. By the construction $C,C'$ have different
singularity types, but their diagrams coincide in all coordinates.

The same applies to the case of $(x^p+y^q+\cdots)(x^p+y^q+\cdots)$: no choice of
coordinates can distinguish this from the pair of branches $(x^p+y^q+\cdots)(x^p-y^q+\cdots)$ (which is
certainly of different type).
\\$\Lleftarrow$ We want to determine the singularity type of $C$ by choosing different coordinates.
For each (\gNnd) branch $C_i$ let $(p_i,q_i)$ be its Puiseux pair (for smooth branches take $(1,\infty)$).
Let $\phi_i$ be a choice of coordinates for which $C_i$ is defined by $x^p_i+y^q_i+\cdots =0$. (The smooth
branches are rectified, i.e. the equation becomes $x=0$).

Let $\Ga_i$ be the \ND~of $(C,0)$ in the coordinates defined by $\phi_i$. Then $\Ga_i$ contains a segment of
slope $p_i/q_i$ (for smooth branches the \ND~is non-commode). Then the collection $\{\phi_i,\Ga_i\}_{i=1}^r$,
for $r-$the number of branches, suffices to determine the singularity type of $(C,0)$.

Indeed, let $C(t)\subset\Sigma_{\{\phi_i,\Ga_i\}_{i=1..r}}$ be a small deformation. By construction,
after the change of coordinates $\phi_i$, the family contains (as a subvariety) the family defined
by $(x^{p_im}\al_1(t)+\cdots+y^{q_im}\al_{m+1}(t))+\cdots=0$. Here inside the terms brackets
correspond to a quasi-homogeneous form, the dots outside correspond to the higher order monomials.
By construction, the quasi-homogeneous form $x^{p_im}\al_1(0)+\cdots+y^{q_im}\al_{m+1}(0)$ is
non-degenerate with $\al_1(0)\neq0$, hence the same holds for $t$ small enough.

Thus, the family $C(t)$ can be decomposed as $\cup C_i(t)$ and the type of each branch is preserved.
(Moreover, all the smooth branches of $C(0)$ stay constant.)
Finally the intersection degrees of the branches are constant (fixed by the set $\{\phi_i,\Ga_i\}$).
\epr

The proof of the proposition gives an upper bound for the number of the needed choices of coordinate system
to restore the singularity type: the number of branches.
This bound can be improved (e.g. in each coordinate system one takes both axes tangent to some branches).
But, e.g. for $r$
singular pairwise non-tangent branches one certainly needs at least $r/2$ coordinate choices.

The proposition also allows to answer positively the question from the introduction: for curves being \gNnd~ or \tNnd~
are properties of singularity types and not only of their representatives.
\bcor\label{Thm tNnd and gNnd is a topological property}
Let $(C,0)\subset(\mC^2,0)$ be a \gNnd~ (or \tNnd) germ of curve. Let $(C',0)\subset(\mC^2,0)$ be a germ of the
same singularity type as $(C,0)$. Then $(C',0)$ is also \gNnd~(or \tNnd).
\ecor
\bpr
$\bullet$ For the \tNnd~ case the statement follows immediately from proposition
\ref{Thm  Criterion on tNnd for curves} (as the conditions are formulated in terms of
 the singularity  types of the branches). Next we consider the other case.

\li
The topological characterization of \gNnd~ curve-germs is well known in the folklore, but we could not find any reference,
except for the preprint \cite{GaLePl07}.

Another way to prove the statement is as follows. By
\cite[Proposition 2.17(2), pg. 287]{GLSbook} the
miniversal equisingular  deformation of any (isolated)
 \Nnd~ plane curve singularity can be realized by monomials not below the \ND, hence consists of \Nnd~ germs.
Let $\mS$ be the singularity type of the \Nnd~germ $(C,0)$, recall that the equisingular stratum $\Si_\mS$
of germs of the singularity type $\mS$ is globally irreducible.
Thus we get: there exists a (Zariski) open dense subset of $\Si_\mS$, whose points correspond to \gNnd~ curve-germs.

On the other hand, for $(C,0)=f^{-1}(0)$, let $\Ga_f$ be the \ND~in the fixed coordinates (so that $f$ is $\Ga_f$
non-degenerate). Consider the stratum of curves that (in {\it some} coordinates)
can be brought to $\Ga_f$ or to a bigger diagram:
\beq
\Si_{\Ga_f}:=\{g\in\mC\{x\}|~~\text{ for some change of coordinates }\phi\circlearrowright(\mC^2,0) ~ ~
\Ga_{\phi^*(g)}\ge\Ga_f\}
\eeq
By definition the stratum is closed, its points correspond either to \gNnd~ curve-germs of type $\mS$
or to higher types adjacent to $\mS$.
The natural morphism $\Si_{\Ga_f}\to\lSi_\mS$ is defined by $g\to (g^{-1}(0),0)$.
By the remark above this morphism is dominant, hence in fact is surjective. Which means:
{\it every} point of $\Si_\mS$ corresponds to a \gNnd~ curve germ.
\epr

\subsection{Germs vs types.}\label{Sec Main Part Germs vs Types}
In this subsection we discuss the questions of \S\ref{Sec Intro Results and Questions} in higher dimensional cases.

For $n\ge3$, being \tNnd~or \gNnd~ is a property of germ representatives (or of analytic
singularity types) but not of the topological types.
We construct equisingular families of surfaces $(V_t,0)\subset(\mC^3,0)$ such that $V_0$ is \Nnd~in
the classical sense but $V_{t\neq0}$ is degenerate (in various senses).

\beR\label{Thm On Bringing to the NND form}
Two observations are useful. Let $(V_f,0)\sset(\mC^3,0)$ with fixed coordinates, such that $f$ is
non-degenerate for $\Ga_f$. So, the coordinates in $\mP^2=Proj(\mC^3)$ are fixed too.
\li  Suppose the projectivized tangent cone $\mP T_{(V_f,0)}\sset\mP^2$ is irreducible, reduced (hence
with isolated singularities). Then its singular locus lies in $[1,0,0]$, $[0,1,0]$, $[0,0,1]\in\mP^2$.
In particular $\mP T_{(V_f,0)}\sset\mP^2$ can have at most three singular points.
\li Let $pt\in Sing(\mP T_{(V,0)})$ and suppose the (plane curve) singularity $(\mP T_{(V,0)},pt)$ is not an \omp.
Let $T=\{l^{p_1}_1\cdots l^{p_k}_k=0\}$ be the tangent cone of $(\mP T_{(V,0)},pt)$, assume $p_1>1$.
Then (as $f$ is \Nnd) $l_1$ coincides with one of the coordinate axes.
\eeR
~\\
The examples below are based on two ideas:
moving several mild singularities of the tangent cone (for the case of \gNnd) or deforming
 one strong singularity of the tangent cone (for the case of \tNnd).
\bex
Consider a super-isolated singularity \cite{Luengo87} $V_0=\{f_p+f_{p+1}=0\}\subset(\mC^3,0)$.
Here $f_{p+1}$ is generic and the projective curve $\{f_p=0\}\subset\mP^2$ has three cusps
(assume $p$ is big enough, e.g. $p\ge6$). According to the remark above,
arrange the coordinates such that the cusps are at  $[1,0,0]$, $[0,1,0]$, $[0,0,1]\in\mP^2$
(corresponding to the coordinate axes $\hat{x}$, $\hat{y}$ and $\hat{z}$ of $\mC^3$).

Note that a $GL(3)$ action which fixes this points is at most a permutation.
To make $V_0$ \Nnd~assume that the tangents to the cusps are oriented along the coordinate axes, e.g.
$f_p=z^{p-3}(zx^2+y^3)+x^{p-3}(xy^2+z^3)+y^{p-3}(yz^2+x^3)$.

Let $V_t$ be the equi-singular family, with the cusps staying at their points $\hat{x},\hat{y},\hat{z}$,
but their tangents changing freely. For example, $f(t)=f_p(t)+f_{p+1}$ with
$f_p(t)=z^{p-3}(z(x-ty)^2+y^3)+x^{p-3}(x(y+tz)^2+z^3)+y^{p-3}(y(z-tx)^2+x^3)$.

Then $V_{t\neq0}$ is \tNnd~but not \gNnd. Indeed, to bring $V_{t\neq0}$ to a \Nnd~ form one should keep the
cusps at the points $\hat{x},\hat{y},\hat{z}$ and at the same time keep their tangents along the axes.
So, only $GL(3)$ transformations are relevant. But, as was noted above, the only $GL(3)$
transformations which keep the cusps at the points $\hat{x},\hat{y},\hat{z}$ are permutations.
\eex

\bex(cf. also \cite[example 5.3]{Altmann91})
Consider the family of surfaces $f_t=f_5(t)+f_6=x^5+z(zx+t y^2)^2+y^5+z^6$.
This is a super-isolated singularity because $Sing(f_5=0)\cap (f_6=0)=\empty$.
The projectivized tangent cone of these surface is the plane
quintic $\{f_5(t)=0\}\subset\mP^2$ with one $A_4$ point at $[0,0,1]\in\mP^2$.
Thus $\mu=68=(5-1)^3+4$ (see \S\ref{Sec Milnor number of tNnd} for the
general formula). The family is equisingular in $t$, e.g. because  each surface $V_t$ is resolved
by one blowup of the origin and the singularity type of exceptional divisor in the surface is independent of $t$.

The singularity $V_{t=0}$ is \Nnd~(by direct check). On the other hand for $t\neq0$ the singularity is not \gNnd.
To show this, we prove that the restriction of $f_t$ to the face $\Ga_{f_t}\cap Span(x^5,y^5,z^5)$
is degenerate for any choice of coordinates. Let $\phi\circlearrowright(\mC^3,0)$ be a local coordinate change.
As we are interested in the face whose monomials correspond to the tangent cone, the non-linear part of $\phi$
is irrelevant. So, assume $\phi\in GL(\mC^3)$ and acts on $\mP T_S=\{f_5=0\}\subset\mP^2$. Thus the
goal is to bring the singularity of this quintic to a \Nnd~form.
But this is impossible for $t\neq0$, because in local coordinates the curve
is defined by $(x+ty^2)^2+x^5+y^5=0$. So, to bring it to a \Nnd~form
we must do a {\it non-linear} transformation $x\to x-ty^2$ on $\mC^2$, which does not
arise from $GL(3)\circlearrowright(\mC^3,0)$.
\eex
\bex{\it Equisingular deformation to essentially \Nd~singularities.}\label{Sec Core Deformation To Essentially Newt. Deg.}
 In the last example all the fibres are \tNnd~ (by corollary \ref{Thm SuperIsolated Sings are tNnd}).
 The next example, in which the generic fibre is essentially \Nd (i.e. not \tNnd), is a simple alteration.
The goal is to construct a \ND~ with the properties: any linear change of coordinates erases
the essential information of the singularity, while any coordinate change whose linear part is
the identity preserves the diagram.

For this, one changes the inclinations of the face on which
the degeneration occurs ($Span(x^5,z^3x^2,zy^4,y^5)$ in the last example) and adds some other faces.
\\
\parbox{13cm}
{Consider the hypersurface $f=x^a+y^b+z^c+z^k(zx+y^2)^2$. (For $(a,b,c,k)=(5,5,6,1)$ one has the previous example.)
Suppose $(a,b,c,k)$ are such that the \ND~ consists of the three faces (cf. the picture): $Conv(x^a,x^2z^{k+2},y^b)$,
$Conv(y^b,x^2z^{k+2},y^4z^k)$ and $Conv(x^2z^{k+2},y^4z^{k},z^c)$. This can be ensured by next convexity conditions:
\li $x^a$ and $z^c$ lie above the plane $Span(x^2z^{k+2},y^4z^k,y^b)$. The equation of this plane
is $\frac{y}{b}+\frac{z}{k}(1-\frac{4}{b})+x(\frac{2}{b}-\frac{1-\frac{4}{b}}{k})=1$, so the conditions are:
\beq
a(\frac{2}{b}-\frac{1-\frac{4}{b}}{k})>1~~\rm{and}~~ \frac{c}{k}(1-\frac{4}{b})>1.
\eeq
}
\begin{picture}(0,0)(-60,0)
\curve(0,0,60,0)\curve(0,0,0,60)\curve(0,0,-40,-40)
\put(52,-10){$y^b$}\put(-35,-48){$x^a$}\put(5,57){$z^c$}
\put(-42,24){$x^2z^{k+2}$}  \put(30,17){$y^4z^k$}

\put(-39.5,-39.5){$\bullet$}\put(48.5,-2){$\bullet$}\put(-12.2,23.1){$\bullet$}\put(19.8,17.9){$\bullet$}
\put(-2.2,55.8){$\bullet$}\put(5,57){$z^c$}

\curve(50,0,-38,-38)\curve(-10,25,-38,-38)\curve(50,0,-10,25)\curve(50,0,22,20)\curve(22,20,-10,25)\curve(22,20,0,58)
\curve(-10,25,0,58)
\end{picture}
\\
\li $z^c$ and $y^4z^k$ lie above the plane $Span(x^a,x^2z^{k+2},y^b)$. The equation of the
plane is $\frac{x}{a}+\frac{y}{b}+\frac{1-\frac{2}{a}}{k+2}z=1$, so the conditions are
\beq
c>\frac{k+2}{1-\frac{2}{a}}~~\text{ and }~~~\frac{k}{1-\frac{4}{b}}>\frac{k+2}{1-\frac{2}{a}}.
\eeq
\li $x^a$ and $y^b$ lie above the plane $Span(z^c,x^2z^{k+2},y^4z^k)$. The equation of the plane
is $\frac{z}{c}+x\frac{1-\frac{k+2}{c}}{2}+y\frac{1-\frac{k}{c}}{4}=1$ giving:
\beq
c>\frac{k+2}{1-\frac{2}{a}}~~\text{ and }~~~\frac{b}{4}(1-\frac{k}{c})>1.
\eeq

Assume further that $b<a<k+4<c$ and also: if $\phi\circlearrowright(\mC^2,0)$ is any locally analytic
transformation whose linear part is identity (i.e. $(x,y,z)\to(x+\phi_x,y+\phi_y,z+\phi_z)$ with $\phi_i\in m^2$)
then $\Ga_f=\Ga_{\phi^*(f)}$. This can be achieved e.g. if for each face the angles with all
the coordinate planes are bounded $\frac{1}{2}<\tan(\al)<2$.

For this it is enough to assume: $z^{c-1}$ lies below the plane $Span(x^a,x^2z^{k+2},y^b)$.

Summarizing, all the restrictions above are implied by the following inequalities:
\beq
b<\frac{k+2}{1-\frac{2}{a}}<\frac{k}{1-\frac{4}{b}}<\min(c,2b),~~~b<a<\min(2b,k+4)<c<k+6.
\eeq
This implies $c=k+5$ and $k>10$. We consider (possibly the simplest case): $f_t=x^{14}+y^{13}+z^{16}+z^{11}(zx+ty^2)^2$.
By direct check this family is equisingular (e.g. $\mu_t=2220=const$, can be calculated using \cite{Singular}).

The generic fibre $f_t^{-1}(0)$ is essentially \Nd by the following proposition.
\bel
Let $f_{t\neq0}$ as above and $g=x^{14}+y^{13}+z^{16}+z^{11}(z^2x^2+y^4)$, \Nnd. Then $f_{t\neq0}$ and $g$ have
the same \ND~ in any coordinate system. But $\mu(g)=2219<\mu(f_t)=2220$, so $f_{t\neq0}$ is not \tNnd.
\eel
\bpr
Let $\phi\circlearrowright(\mC^2,0)$ be a locally analytic change of coordinates whose linear part is identity.
By the construction it preserves the \ND. Therefore it is enough to consider only {\it linear} coordinate changes.
But then only the monomials $x^a,y^b$ are relevant and their coefficients are the same in both cases.
\epr
\beR
\li The families in the examples above are not just $\mu-$constant but $\mu^*-$constant (by direct computation).
In fact, one can show that in these families all the polar multiplicities are preserved,
i.e. the surfaces are c-cosecant in the sense of \cite{Teissier77}.
\li In the examples above the singular types admit some \Nnd~representatives, but the generic
representatives are not \gNnd. More precisely, consider the equisingular stratum in the space of miniversal deformation.
The locus of \gNnd~surface-germs is of positive codimension.
\eeR
\eex

\section{\DNnd~singularities}\label{Sec DNnd singularities}
As was shown above,  the \tNnd~ germs form quite a broad class, difficult to work with. One could consider
other intermediate classes of singularities (between \gNnd~and \tNnd). We study here the minimal generalization: \dNnd~germs.
These are higher dimensional analogs of example \ref{Ex Intro Collection Newton Diagrams}.

For a fixed diagram $\Ga$ define the subset $\De\subset\Ga$ of faces "far from the $\hat{x}_n$ axis" as follows.
Let $\hat{x}_n$ be the unit normal to the hyperplane $Span(x_1,\cdots, x_{n-1})$.
For each top-dimensional face $\si\subset\Ga$ let $Span(\si)$ be the supporting hyperplane.
 Let $\hat{v}_\si$ be the unit normal to $Span(\si)$,
oriented such that $(\hat{x}_n,\hat{v}_\si)>0$.

We define:
\beq
\De:=\overline{\cup \Big\{\bM\si\subset\Ga\\\text{\tinyT top-dimensional}\eM|~~(\hat{x}_n,\hat{v}_\si)\ge\frac{1}{\sqrt{n}}\Big\}}.
\eeq
Since the union is over the top dimensional faces, $\De$ can be empty.
Here the topological closure is needed to add the relevant non-top-dimensional faces.
\bex
In the case of plane curves, with $\Ga\subset\mR^2_{\ge0}$, one gets that  $\De$ is the union of all the faces
whose angle with the $\hat{x}_1$ is not bigger than $\pi/4$.

For the \ND~of the surface $x_1x_2x_3+x^p_1+x^q_2+x^r_3$ with $p,q,r\ge3$ one has: $\De=Conv(x_1x_2x_3,x^p_1,x^q_2)$.

For the \ND~of the surface $x^N_3+x^p_2+x^p_3$ with $N>p$ one has: $\De=\empty$.
\eex
 In the sequel we assume that the projectivized tangent cone $\mP T_{(V_f,0)}$ has only isolated singularities,
in particular it is reduced.
\bed\label{Def dNnd}
The isolated hypersurface singularity $(V_f,0)\subset(\mC^n,0)$ is \dNnd~if for each singular
point $x\in Sing(\mP T_{(V_f,0)})$
there exists a coordinate-change $\phi\circlearrowright(\mC^n,0)$ such that $\phi(x)=[0,\cdots ,0,1]$
and the restriction
$\phi^*(f)|_{\overline{\Ga\smin\De}}$ is non-degenerate.
\eed
\bel
(\bet generalized\\Newton-non-degenerate\\with $Sing(\mP T_{(V_f,0)})$ isolated\eet)$\Rrightarrow$(\bet directionally\\Newton-non-degenerate\eet)
$\Rrightarrow$(\bet topologically\\Newton-non-degenerate\eet).
\eel
\bpr
The first implication is obtained as follows. Suppose $f$ is \Nnd~in some fixed coordinates
and $pt\in Sing(\mP T_{(V_f,0)})\sset\mP^{n-1}$. Let $p=mult(V_f,0)$. By the non-degeneracy on the
face $\Ga_f\cap Span(x^p_1,\cdots,x^p_n)$ one gets that  $pt\not\in(\mC^*)^{n-1}$, so by a permutation
of coordinates we can assume: $pt=[0,*,\cdots,*]\in\mP^{n-1}$. Now consider the restriction of $f$ to
the boundary $\di\Big(\Ga_f\cap Span(x^p_1,\cdots,x^p_n)\Big)$. Checking the non-degeneracy on
the top-dimensional components one has (after a permutation) $pt=[0,0,*,\cdots,*]\in\mP^{n-1}$. After several similar steps we get that  $pt=[0,0,\cdots,0,1]$.
\\
\\
For the second implication suppose $(V_f,0)$ is \dNnd. Consider the corresponding stratum
\beq
\Si_f=\overline{\{g|~\Ga_g=\Ga_f ~ \text{ in all coordinate systems}\}}.
\eeq
\li First we prove that any family $(V(t),0)\subset\Si_f$ ~ for $t$ small enough consists of \dNnd~hypersurfaces.
Recall that the tangent cone is constant in this family (cf. proposition \ref{Thm Tangent Cone of tNnd is fixed}).
For any point $pt\in Sing(\mP T_{V(0)})=Sing(\mP T_{V(t)})$ choose coordinates with $pt=[0,..,0,1]$
 and $f(0)|_{\overline{\Ga\smin\De}}$ non-degenerate.

Note that $\Ga_{f(t)}=\Ga_{f(0)}$ and for the fixed \ND~the non-degeneracy is an open property in Zariski topology.
So for each singular point $pt\in Sing(\mP T_{V(0)})$ we have a bound: if $|t|<\ep(pt)$
then $f(t)|_{\overline{\Ga\smin\De}}$ is non-degenerate.
Thus, as $\mP T_{V(0)}$ has isolated singularities only, the needed bound is found by combining a finite
number of inequalities.

\li
Now observe that such a family $(V(t),0)\subset\Si_f$, for $|t|<\ep$,
consisting of \dNnd~germs, is equi-resolvable in the following sense.
Blowup the origin, let $\tilde{V}(t)$ be the strict transform. By construction it has a finite
number of singular points $pt\in Sing(\mP T_{V(0),0})\sset\mP^{n-1}$, their positions are independent of $t$. At each
such point the germ $(\tilde{V}(t),pt_\al)$ is \Nnd~ in some local coordinates, the choice of
the coordinates is independent of $t$, and in the chosen coordinates the \ND~ $\Ga_{(\tilde{V}(t),pt_\al)}$
is independent of $t$. Now we can use various nice properties of families of \Nnd~ hypersurfaces \cite{Oka79}.

Apply the toric modification to the multigerm $\coprod_\al (Bl(\mC^n),pt_\al)$,
corresponding to the diagrams $\Ga_{(\tilde{V}(t),pt_\al)}=\Ga_{(\tilde{V}(0),pt_\al)}$.
Let $(Y,E)\to(\mC^n,0)$ be the resulting total modification of the original ambient space.
Now the strict transform $\widetilde{\widetilde{V}}(t)$ of $(V(t),0)$ is a family of {\it smooth} hypersurfaces,
intersecting transversally the exceptional divisor $E$. Moreover, for each irreducible
component $E_i\sset E$ the intersections $C(t)=\widetilde{\widetilde{V}}(t)\cap E_i$
form a family of embedded varieties $(E_i,C(t))$, which is locally trivial over the base $(\mC^1_t,0)$.
In particular, for any small $t$, one has embedded diffeomorphism $(E_i,C(t))\approx (E_i,C(0))$.

\li
Now we prove that such a family is $\mu$-const.
The Milnor number is determined by the monodromy zeta function, hence it is enough to show that the zeta
function is constant. And the later is immediate by the theorem of \cite{A'Campo75}, we recall it below
\S\ref{Sec Zeta function for tNnd}, equation (\ref{Eq Monodromy ACampo}).
Indeed, by the above discussion, for each $m$ we have that $\chi(S_m)$ is independent of $t$.
\\
Finally, observe that $\mu$-constant implies equisingularity, for $n=3$ note that the family can always
be chosen linear.
\epr

 The three non-degeneracy classes are distinct as is seen already in the case of curves:
\bprop
The plane curve germ $(C,0)\subset(\mC^2,0)$ is \dNnd~iff in its tangential decomposition $(C,0)=\cup (C_\al ,0)$ each $(C_\al ,0)$ is \gNnd.
In particular, being \dNnd~is a property of the singularity type. Namely, if $(C,0)$ is \dNnd~then any other germ of the same type is \dNnd.
\eprop
In particular, for curves being \dNnd~places no conditions on the tangent cone.
In higher dimensions the tangent cone is more restricted.
\bprop\label{Thm Tangent Cone for dNNd}
Let $(V_f,0)\subset(\mC^n,0)$ be \dNnd.
\bee
\item \label{111}
Every singular point of $\mP T_{(V_f,0)}\subset\mP^{n-1}$ can be brought to a \Nnd~ form by a linear transformation (the corresponding \ND~can be non-commode).
\item Let $f=f_p+f_{p+k}$ with $Sing(\mP T_{V_f})$ isolated and $f_{p+k}$ generic with respect to $f_p$
(i.e. $Sing(V_{f_p})\cap V_{f_{p+k}}=\empty$).
If $f_p$ (i.e. $\mP T_{(V_f,0)}$) satisfies the condition (\ref{111}) above
then $(V_f,0)$ is \dNnd.
\eee\eprop
\bpr
\bee
\item Let $pt\in Sing(\mP T_{V_f,0})$, choose the coordinates such that $pt=[0,\cdots,0,1]$
and $f|_{\overline{\Ga\smin\De}}$ is non-degenerate. Let the Taylor expansion be $f=f_p+\cdots$. Note that the \ND~
of the singularity $(f^{-1}_p(0),0)\subset(\mP^{n-1},0)$ is obtained from $\Ga_{(V_f,0)}$ by intersecting
with the hyperplane $\sum x_i=p$ and projecting to the coordinate hyperplane $Span(x_1,\cdots,x_{n-1})$. Thus,
in the fixed coordinates, $f_p$ is non-degenerate for its \ND. Finally, observe that $f_p$ was brought to
this form by a $GL(n)$ transformation, because for any locally
analytic change of coordinates in $(\mC^n,0)$ only the linear part acts nontrivially on the tangent cone.
\item
Let $pt\in Sing(\mP T_{V_f,0})\sset\mP^{n-1}$. As the transformation, which brings $(\mP T_{(V_f,0)},pt)$ to a
non-degenerate form, is linear, it can be lifted to a  linear transformation of $(\mC^n,0)$.
Then, the restriction of the transformed $f$ to the relevant boundary components of $\Ga_f\cap Span(x^p_1,\cdots,x^p_n)$
is non-degenerate. Hence, by the genericity of $f_{p+k}$, the restriction $f|_{\overline{\Ga\smin\De}}$ is non-degenerate.
\epr
\eee
\subsection{The directional approximations}\label{Sec tNnd directional approximation}
We want to generalize the directional approximations of curves (introduced in example
\ref{Ex Intro Collection Newton Diagrams}).
The reader is invited to review the definition of "faces far from the $\hat{x}_n$ axis".
\\\parbox{14cm}
{
\bed\label{Def Directional Approximation}
Let $(V,0)\subset(\mC^n,0)$ be \dNnd, in particular $\mP T_{(V,0)}\subset\mP^{n-1}$ has isolated singularities only.
We say that $(V,0)$ has a \Nnd~directional approximation $(V_1,0)$, $\ldots$, $(V_k,0)$ if for any
point $pt_\al\in Sing(\mP T_{(V,0)})$
there exists a coordinate system  $\phi_\al \circlearrowright (\mC^n,0)$, such that:
\ls $\phi_\al(pt_\al)=[0,\cdots,0,1]\in\mP^{n-1}$ and the subset $\De$ of faces "far from the $\hat{x}_n$ axis" lies inside the hyperplane $\sum x_i=mult(V,0)$
\ls $\De\cap Span(\hat{x}_1,\cdots,\hat{x}_{n-1})=Conv(x^p_1,\cdots,x^p_{n-1})$ for $p=mult(V,0)$
\ls $f|_{\overline{\Ga_{\phi^*_\al(f)}\smin\De}}$ is \Nnd, $\Ga_{(V_\al,0)}=\Ga_{(\phi^*_\al(f))}$
and $(V_\al,0)$ is \Nnd.
\eed
}
\begin{picture}(0,0)(-40,20)
\put(0,0){\vector(0,1){70}}\put(0,0){\vector(2,-1){50}}\put(0,0){\vector(-1,-1){30}}
\put(-22.5,-22.5){$\bullet$}\put(27.5,-17.5){$\bullet$}\put(-12.5,-7.5){$\bullet$}\put(17.5,-2.5){$\bullet$}

\put(-15,65){$\hat{x}_n$}\put(-25,-35){$\hat{x}_1$}\put(35,-30){$\hat{x}_{n-1}$}

\thicklines
\curve(-20,-20,-10,-5)\curve(30,-15,20,0)\curve(30,-15,-20,-20)
\curve(-10,-5,-2,3)\curve(7,4,-2,3)\curve(7,4,20,0)
\curve(0,60,-10,-5)\curve(0,60,20,0)\curve(0,60,7,4)\curve(0,60,-2,3)
\thinlines
\put(-27,-20){$\tinyA p$}\put(33,-13){$\tinyA p$}
\put(-30,-12){$\De$}\put(-22,-10){\vector(1,0){20}}
\end{picture}
\\
Note that the existence of a
\Nnd~directional approximation, in general, imposes additional restriction for $f$:
\bex
$\bullet$ Any \dNnd~plane curve singularity has a \Nnd~directional approximation (cf. example \ref{Ex Intro Collection Newton Diagrams}).
\li On the other hand, this is not the case if $n\geq 3$. For such $n$, a  sufficient condition is the following. Let $(V,0)\subset(\mC^n,0)$ be \dNnd,
with $Sing(\mP T_{(V,0)})=\{pt_1,\cdots,pt_k\}\sset\mP^{n-1}$. Suppose
the \ND~of  any germ $(\mP T_{(V,0)},pt_\al)\sset(\mC^{n-1},z_i)$ remains commode after any $GL(n-1)$
change of coordinates (preserving $pt_\al$).
Then there exists a \Nnd~directional approximation.
\eex
Similarly to the case of curves, the germs $(V_i,0)$ are not defined uniquely. But any two representatives
(associated with a fixed point $p_\alpha$) have
the same singularity type, the same \ND~ and can be joined by a $\mu-$constant family.
\section{Some numerical singularity invariants}\label{Sec Invariants}
Here we generalize some classical theorems to \dNnd~germs which admit  \Nnd~directional approximation
(cf. definition \ref{Def Directional Approximation}). The proofs of the statements are
 in \S \ref{Sec Singularity Invariants The Case of Curves} and
 \S\ref{Sec Singularity Invariants The Case of Hypersurfaces}.
\subsection{Kouchnirenko's formula for the Milnor number}\label{Sec Milnor number of tNnd}
Recall, that for a \Nnd~germ the Milnor number is determined by the \ND~by the classical Kouchnirenko formula:
\beq\label{Eq Kouchnirenko Formula}
\mu(V,0)=\sum_{j=0}^n(-1)^{n-j}j!Vol_j(\Ga_{(V,0)}),
\eeq
where  $Vol_j(\Ga_{(V,0)})$ is the sum of $j-$dimensional volumes of the intersections of $\Ga$ with
the $j-$dimensional coordinate hyperplanes. The volume of unit cube is 1 in any dimension.

\bthe\label{Thm Milnor Number Kouchnirenko} Let $(V,0)\subset(\mC^n,0)$ be an isolated \dNnd~hypersurface singularity with multiplicity $p$.
Let $Sing(\mP T_{(V,0)})=\{pt_1, \ldots, pt_k\}$ and let $(V_1,0),\cdots, (V_k,0)$ be its \Nnd~directional approximations
 (cf. \S\ref{Sec tNnd directional approximation}).
Then Kouchnirenko's formula holds in the following form:
\beq\label{uj}
\mu(V,0)=\sum_{\al=1}^k\mu(V_\al,0)-(k-1)(p-1)^n,\text{ where }~\mu(V_\al,0)\text{ is computed by the
formula }(\ref{Eq Kouchnirenko Formula}).
\eeq
\ethe
\bex\label{Ex Milnor number for dNnd}
$\bullet$ For plane curves this proposition is just a reformulation of
tangential decomposition, see \S\ref{Sec Singularity Invariants The Case of Curves}.
Let $(C,0)=\cupl^k_{\al=1}(C_\al,0)$ be the tangential decomposition, let $\mu_\al,p_\al$ be the Milnor
number and the multiplicity of $C_\al$, so that $p=\sum p_\al$.
Then, for the directional approximation one has: $\mu(C^{(\al)})=\mu_\al+(p-p_\al)(p+p_\al-2)$.
Thus, by (\ref{uj}), the total Milnor number is:  $\mu=\sum \mu_\al+p^2+1-k-\sum p^2_\al$.
\li Let $(V_f,0)\sset(\mC^n,0)$ for $f=f_p+f_{p+1}+\cdots\in\mC\{x\}$ be \dNnd~with the directional
approximation $\{(V_\al,0)\}$  corresponding to the singular points $\{pt_\al\in\mP T_{(V,0)}\}$.
Assume each $V_\al$ has the singularity type of the germ
 $\{f_p^{(\al)}(x_1,\cdots,x_n)+x^{p+q_\al}_n=0\}\subset \mC^n$, where the hypersurface
 $\{f_p^{(\al)}(x_1,\cdots,x_n)=0\}\sset\mP^{n-1}$ is singular at the point $[0,\cdots,0,1]\in\mP^{n-1}$ only.
 For example, \dNnd~ singularities of Yomdin type ($f_p+f_{p+k}$) are of this kind.
By the theorem the computation of $\mu(V_f,0)$ reduces to $\mu(V_\al,0)$.

Blowup $\mC^n$ at the origin, let $\tilde{V}_\al$ be the strict transform, then (see \cite{Melle00} or
the proof of this theorem):
\beq
\mu(V_\al,0)=(p-1)^n+\mu(\mP T_{(V_\al,0)},pt_\al)+\mu(\tilde{V}_\al).
\eeq
Note that $\tilde{V}_\al$ is the suspension
(in some local coordinates its equation is of the form: $g(x_2,\ldots,x_{n-1})+x^{q_\al}_n=0$) thus
$\mu(\tilde{V}_\al)=\mu(\mP T_{(V_\al,0)},pt_\al)(q_\al-1)$ and we obtain:
\beq
\mu(V,0)=(p-1)^n+\sum q_\al\mu(\mP T_{(V,0)},pt_\al).
\eeq
If all the $q_\al$ are equal this reproduces the result of \cite{Luengo-Melle95} (in the \dNnd~ case).
\li It is instructive to compute some case explicitly. For example, consider $f=f_4+f_6+f_8$, where
\beq
f_4=xyz(x+y+z), ~ ~ f_6=(x^2y^2+y^2z^2+x^2z^2)(x+2y+3z)^2, ~ ~ f_8=x^8+y^8+z^8
\eeq
Then the curve $f^{-1}_4(0)\sset\mP^2$ (four lines) has nodes at
the points $[1,0,0]$, $[0,1,0]$, $[0,0,1]$, $[1,1,0]$, $[1,0,1]$, $[0,1,1]$.
The directional approximations at  $[1,1,0]$, $[1,0,1]$, $[0,1,1]$ are of the singularity
type $x^4+y^4+x^2z^2+y^2z^2+z^6$. So $q_\al=2$ and $\mP T_{(V_\al,0)}$ has only one singular point: the node,
so $\mu(\mP T_{(V_\al,0)})=1$.
The directional approximations at $[1,0,0]$, $[0,1,0]$, $[0,0,1]$ have the singularity type of
 $x^4+y^4+x^2z^2+y^2z^2+z^8$. So $q_\al=4$, while $\mP T_{(V_\al,0)}$ has again only one node,
 $\mu(\mP T_{(V_\al,0)})=1$.

In total one has:
\beq
\mu=3^3+(2\times 1+2\times 1+2\times 1)+(4\times 1+4\times 1+4\times 1)=45
\eeq
which can be of course recomputed directly, e.g. by \cite{Singular}.
\eex
\subsection{Zeta function of monodromy}\label{Sec Zeta function for tNnd}
Recall the basic result of \cite{A'Campo75} (see also \cite[II.3.12]{AGLV}).
\\
\parbox{14cm}
{ Given an isolated hypersurface singularity, construct its good embedded resolution (see the diagram): $\tV$ is
smooth, $E$ consists of smooth components
and $\tV\cup E$ is a normal crossing divisor. Let $\pi^{-1}(0)=\sum m_iE_i$, i.e. $E_i$ is an irreducible
component of $E$ of multiplicity $m_i$. Define $S_m:=\{x\in E^{reg}:~mult(E,x)=m\}$ (where $E^{reg}$ is the regular part of $E$).
}\hspace{0.3cm}
$\bM
(\tV,\tV\cap E)&\subset(Y,E)\\\downarrow&\downarrow\pi\\
(V,0)&\subset(\mC^n,0)
\eM$
\\
 Then
\beq\label{Eq Monodromy ACampo}
\ze_{(V,0)}(z)=\prod_{m\ge1}(1-z^m)^{\chi(S_m)}
\eeq
where $\chi$ is the Euler characteristic.

The product structure of this formula is the basic reason for the possibility to determine the zeta function
by the geometry of the \ND.
\bprop\label{Thm Zeta Function for tNnd}
Let $(V,0)\subset(\mC^n,0)$ be a \dNnd~germ, whose projectivized tangent cone has isolated singularities:
$Sing(\mP T_V)=\{pt_1,\cdots,pt_k\}$, let $(V_1,0),\ldots, (V_k,0)$ be the corresponding directional approximations.
Then A'Campo's formula can be written in the form
\beq
\ze_{(V,0)}(z)=\frac{\prod\limits_{\al=1}^k\ze_{(V_\al,0)}(z)}{(1-z^p)^{(k-1)(n-\chi_{p,n-1})}}
\eeq
where $\ze_{(V_\al,0)}(z)$ is the classical zeta-function of the (\Nnd) hypersurface-germ,
$p=mult(V,0)$ and \mbox{$\chi_{p,n-1}=\chi(V_{p,n-1})=\frac{(1-p)^{n}-1}{p}+n$} is the topological Euler characteristic of
an arbitrary smooth hypersurface $V_{p,n-1}\subset\mP^{n-1}$ of degree $p$.
\eprop
For a related result, see \cite{Gusein-Zade Luengo Melle-Hernandez97}.

Using the last proposition it is immediate to generalize Varchenko's formula for the zeta function in terms of the \ND.
Recall (\cite{Varchenko76},  \cite[II.3.12]{AGLV}) that for an isolated \Nnd~hypersurface singularity $(V,0)\subset(\mC^n,0)$ the zeta
function of the monodromy can be written in the form:
\beq
\ze_{(V,0)}(z)=\prod_{l=1}^n(\ze^l(z))^{(-1)^{l-1}},
\eeq
where $\{\ze^l(z)\}_l$ are some polynomials completely determined by the geometry of $l-$dimensional faces of the \ND.
\bcor\label{Thm Zeta Function for tNnd in terms of ND}
Under the assumptions of proposition \ref{Thm Zeta Function for tNnd} Varchenko's formula is valid in the following form:
\beq
\ze_{(V,0)}(z)=\frac{1}{(1-z^p)^{(k-1)(n-\chi_{p,n-1})}}\prod_{pt_\al\in Sing(\mP T_{(V,0)})}
\prod_{l=1}^n(\ze^l_{(V_\al,0)}(z))^{(-1)^{l-1}}
\eeq
where, \mbox{$\chi_{p,n-1}=\chi(V_{p,n-1})$} is the topological Euler characteristic of an arbitrary smooth hypersurface
$V_{p,n-1}\subset\mP^{n-1}$ of degree $p$, and $\ze^l_{(V_\al,0)}(z)$ are the standard Varchenko
polynomials (for the \Nnd~ singularities $V_\al$).
\ecor
\bex\label{Ex Zeta Function For dNnd} $\bullet$
Let $(C,0)\subset(\mC^2,0)$ be the union of $k$ non-tangent branches of types $x^{p_\alpha}+y^{q_\alpha}$ for $p_\alpha<q_\alpha$,
 $(p_\alpha,q_\alpha)=1$. Let $p=mult(C,0)=\sum p_\alpha$.
The proposition \ref{Thm Zeta Function for tNnd} gives: $\ze=\frac{\prod\ze_\alpha}{(1-z^{p})^{(k-1)(2-p)}}$
where $\ze_\alpha=\frac{1-z^{p+q_\alpha-p_\alpha}}{(1-z^{p})^{p-p_\alpha-1}(1-z^{p_\alpha(p+q_\alpha-p_\alpha)})}$.
Altogether one has: $\ze=(1-z^p)^{2-k}\prod_\alpha\frac{1-z^{p+q_\alpha-p_\alpha}}{1-z^{p_\alpha(p+q_\alpha-p_\alpha)}}$.

\li As in the example \ref{Ex Milnor number for dNnd}, consider a \dNnd~ hypersurface
\beq
(V_f,0)\sset(\mC^n,0) \text{ for } f=f_p+f_{p+1}+\cdots\in\mC\{x\}
\eeq
with the directional
approximation $\{(V_\al,0)\}$  corresponding to the singular points $\{pt_\al\in\mP T_{(V,0)}\}$.
Assume each $V_\al$ has the singularity type of the germ
 $\{\sum_{i<n} x^p_i+x_n^{p-p_\al}\sum_{i<n} x^{p_\al}_i +x^{p+q_\al}_n=0\}\sset\mC^n$, i.e.
 its projectivized tangent cone has only one singular point: an \omp~ of multiplicity $p_\al$.
 For simplicity assume $(p_\al,q_\al)=1$.

Then, for $n=3$,  an immediate application of Varchenko's formula gives:
\beq
\ze(V_\al,0)=\frac{(1-z^p)^{p^2-p^2_\al}(1-z^{p_\al(q_\al+p)})^{p_\al}(1-z^{(q_\al+p)})(1-z^p)^2}
{(1-z^p)^{2(p-p_\al)}(1-z^{p_\al (q_\al+p)})^2(1-z^p)^p}=
(1-z^p)^{p^2-3p+3-(p_\al-1)^2}(1-z^{p_\al (q_\al+p)})^{p_\al-2}(1-z^{(q_\al+p)})
\eeq
Substituting this into the formula of proposition \ref{Thm Zeta Function for tNnd} one has:
\beq
\ze_{(V_f,0)}=(1-z^p)^{p^2-3p+3-\suml_\al(p_\al-1)^2}\prod_\al(1-z^{(q_\al+p)})(1-z^{p_\al (q_\al+p)})^{p_\al-2}
\eeq
As a consistency check we can get the Milnor number, as the degree of the polynomial $\frac{\ze_{(V_f,0)}}{1-z}$.
From the last formula the degree is:
\beq
p(p^2-3p+3)-1-\sum p(p_\al-1)^2+\sum_\al\Bigg( (q_\al+p)+p_\al (q_\al+p)(p_\al-2)\Bigg)=(p-1)^3+\sum q_\al(p_\al-1)^2
\eeq
recovering the result of example \ref{Ex Milnor number for dNnd} for our case.
\eex
\subsection{Order of determinacy}
Suppose a \dNnd~ hypersurface germ $\{f=0\}=(V,0)\subset(\mC^n,0)$ has a \Nnd~directional
approximation $(V_1,0),\cdots, (V_k,0)$. For each $(V_\al,0)$ let $o.d.(V_\al,0)$ be the (contact, topological)
order of determinacy \cite[I.2.2]{GLSbook}. It is easily read from the diagram of $(V_\al,0)$.
\bprop
The order of determinacy of $(V,0)$ is $\underset{\al}{max}(o.d.(V_\al,0))$.
\eprop
\bpr
Let $q=\underset{\al}{max}(o.d.(V_\al,0))$. Suppose $jet_q(f)=jet_q(g)$, then $\Ga_f=\Ga_g$ in any coordinates.
Hence $g$ is \tNnd~with respect to its collection of the diagrams. So, the singularity
types of $f,g$ coincide.

On the other hand the order of determinacy of $(V,0)$ is certainly at least $\underset{\al}{max}(o.d.(V_\al,0))$.
\epr
\subsection{Proofs. The case of curves}\label{Sec Singularity Invariants The Case of Curves}
For curves the proofs are especially simple (and they will be related with some other invariants as well, e.g. with the modality).
\subsubsection{Proof of Proposition \ref{Thm Milnor Number Kouchnirenko}}
It is based on the formulas for the $\de$ invariant:
\li $\mu=2\de-r+1$, here $r$ is the number of branches
\li $\de=\suml_\al\de(C_\al)+\suml_{\al<\be}<C_\al,C_\be>$ (for the tangential decomposition $C=\cupl_{\al=1}^k C_\al $).
\\
Using these formulas one gets:
\beq
\mu(C)=\suml_\al(\mu(C_\al)+r_\al-1)+2\suml_{\al<\be}<C_\al,C_\be>-\suml_\al r_\al+1=
\suml_\al\mu(C_\al)+1-k+\suml_{\al\neq \be}p_\al p_\be
\eeq
(for $p_\al=mult(C_\al)$). Assume that the curve-germ has \Nnd~directional approximations $\{C^{(\al)}\}$
(i.e. each $C_\al$ is \gNnd). Then the result follows from the observation:
$\mu(C^{(\al)})=\mu(C_\al)+p^2-p^2_\al+2p_\al-2p$.

\subsubsection{Proof of Proposition \ref{Thm Zeta Function for tNnd}}
\label{Sec Zeta func. for curves tNnd}
Let $C=\cupl_{\al=1}^k C_\al $ be the tangential decomposition, then the resolution tree consists of $pE$
and $k$ chains corresponding to $\{C_\al \}$. Let $\{C^{(\al)}\}$ be the directional approximations
(cf. example \ref{Ex Intro Collection Newton Diagrams}).
The embedded resolution graph of $(C,0)$ is obtained as the union of those of $\{C^{(\al)}\}$, they
all have the same component $pE$. So, the product in equation \ref{Eq Monodromy ACampo}
splits into parts corresponding to $\{C^{(\al)}\}$ and one should divide by the over-counting:
\beq
\ze_C(z)=\frac{\prod_\al\ze_{C^{(\al)}}(z)}{(1-z^p)^{(k-1)(2-p)}}
\eeq
Note that in the case of curves this result is valid without assumption that each $C_\al $ is \gNnd.
\subsubsection{} Corollary \ref{Thm Zeta Function for tNnd in terms of ND} can be stated very explicitly in
the \dNnd~case.
\\
\parbox{13cm}
{Assume each $C_\al$ is \gNnd. To write Varchenko's formula introduce the parameters of the \ND~of $\Ga_{C^{(\al)}}$.
For each edge $l_j$ of the diagram (except for $\overline{(0,p),(p-p_\al,p_\al)}$) let $a_j x+b_j y=c_j$
be the equation of the line it spans. Here $a_j,b_j,c_j\in\mN$ and $(a_j,b_j)=1$
(so the coefficient are fixed uniquely). Let $|l_j|$ be the number of integral points minus one, on
the edge $j$. Then the formula of corollary \ref{Thm Zeta Function for tNnd in terms of ND} reads:
}
\begin{picture}(0,0)(-20,20)
\put(0,0){\vector(1,0){90}}\put(0,0){\vector(0,1){50}}
\put(18,18){$\bullet$}\put(-2.5,37.5){$\bullet$}\put(37.5,7.5){$\bullet$}\put(77.5,-2.5){$\bullet$}

\curve(0,40,20,20)\curve(20,20,40,10)\curve(40,10,80,0)
\put(-8,38){$\tinyA p$}\put(-8,18){$\tinyA p_\al$}\put(13,-7){$\tinyA p-p_\al$}\put(78,-7){$\tinyA q$}
  \put(70,40){$ \Ga_{C^{(\al)}}$}      \put(10,40){ignore}\put(20,38){\vector(-1,-1){10}}

\put(48,10){$\tinyA a_{j_\al} x+b_{j_\al} y=c_{j_\al}$}
\end{picture}
\\
\beq
\ze_C(z)=\frac{1}{(1-z^p)^{k-2}}\prod_{\al=1}^k\frac{1-z^{q_\al}}{\prod_{\al_\al}(1-z^{c_{\al_\al}})^{|l_{\al_\al}|}}
\eeq
Here $\al_i$ runs over the edges of $\Ga_{C^{(i)}}$, each time omitting the edge $\overline{(0,p),(p-p_i,p_i)}$.

\subsection{Right modality for functions of two variables}\label{Sec.Right Modality in two variables}~\\
\parbox{13cm}
{For the \Nnd~singularities the (right) modality can be calculated  as the number of integral points $(x,y)$
not above the \ND, satisfying $x,y\ge2$ (cf. \cite{Kouchnirenko76}). This generalizes naturally to \dNnd~case.
}
\begin{picture}(0,0)(-20,10)
\put(0,0){\vector(1,0){90}}\put(0,0){\vector(0,1){50}}
\put(18,18){$\bullet$}\put(-2.5,37.5){$\bullet$}\put(37.5,7.5){$\bullet$}\put(77.5,-2.5){$\bullet$}
\multiput(18,18)(3,-3){8}{$\cdot$}
\curve(0,40,20,20)\curve(20,20,40,10)\curve(40,10,80,0)
\put(-8,38){$\tinyA p$}\put(-8,18){$\tinyA p_\al$}\put(13,-7){$\tinyA p-p_\al$}  \put(30,25){$ \Ga_{C^{(\al)}}$}

\put(60,15){$\De_\al$}  \put(60,13){\vector(-1,-1){10}}
\end{picture}
\\
Let $C=\cupl_\al C_\al $ be the tangential decomposition, assume each $C_\al $ is \gNnd.
Let $\Ga_{C^{(\al)}}$ be the \ND~of the corresponding directional approximation.
The {\it closed} set $\Ga_{C^{(\al)}}^-$  of points not above the \ND~ decomposes naturally into the triangle $x+y\le p$ for $p=mult(C)$ and
the remaining polygon $\De_\al$. Note that $\De_\al$ is `half--open': if $(x,y)\in \De_\al$ then $x+y>p$.
\bprop
$\bullet$ In the notations above, the right modality of $(C,0)$ is $\bin{p-2}{2}+\suml_\al\sharp|\De_\al\cap\mZ^2_{\ge2}|$, where $\sharp$ counts the number of lattice points.
\li Alternatively, if $(\mu_\al,p_\al,modality_\al)$ are invariants of $C_\al$ then
modality(C)=$\sum modality_\al+\frac{p^2}{2}-\sum\frac{p^2_\al}{2}+3-3k$.
\eprop
In the last formula the convention is: the right modality of $A_0$ (i.e. of a smooth branch) is 1.
\\
\bpr
By definition, the right modality is $\mu-\tau^{es}$, where $\tau^{es}$ is the codimension of the
equi-singular stratum in the space
of miniversal deformation. The additivity of $\mu$ for the diagrams of directional approximation follows from
proposition \ref{Thm Milnor Number Kouchnirenko}.

For $\tau^{es}$ use the classical formula (in terms of the multiplicities of the strict transform on the resolution),
cf. \cite[pg. 373, eq. 2.8.36]{GLSbook}
\beq
\tau^{es}=\sum\bin{m_q+1}{2}-1-\sharp\Big({\ber free\\points\eer}\Big).
\eeq
Here the sum is over the infinitely near points of $(C,0)$ arising during the resolution. An infinitely
near point of $(C,0)$ is free if it lies on at most one exceptional divisor.
In particular the initial singular point of $(C,0)$ is free. For a smooth branch the convention is $\tau^{es}(A_0)=-1$.

Finally, as in the proof for zeta function, relate the resolution tree of $(C,0)$ to the trees
of its directional approximations $(C^{(\al)},0)$.

The second claim follows similarly using: $\tau^{es}=\sum_\al\tau^{es}(C_\al)+\bin{p+1}{2}-\sum\bin{p_\al+1}{2}+2k-2$.
Note the convention in this formula: $\tau^{es}(A_0)=-1$.
\epr
\bex
Let $(C,0)$ be the union of $\{A_{n_i}\}_{i=1}^k$ with tangent lines pairwise distinct. Then for
each directional approximation $\sharp|\De_\al\cap\mZ^2_{\ge2}|=0$. Therefore, the
right modality is $\bin{2k-2}{2}$.
\eex
\subsection{Proofs. The case of hypersurfaces}\label{Sec Singularity Invariants The Case of Hypersurfaces}
\subsubsection{Proof of Theorem \ref{Thm Milnor Number Kouchnirenko}}
As the Milnor number is determined by the zeta function of monodromy
(since $\deg{\zeta}=1+(-1)^{n-1}\mu$)), the formula can be immediately
obtained from the proposition \ref{Thm Zeta Function for tNnd}. We give also a direct proof.

The derivation of the formula is based on the following result \cite[Theorem 1]{Melle00}.
For the hypersurface germ $(V,0)\subset(\mC^n,0)$, let $\mP T_{(V,0)}\subset\mP^{n-1}$  be
the projectivization of its tangent cone and $\tV\to V$ the strict transform under the blow-up of the origin.
Assume, both $\tV$ and $\mP T_{(V,0)}$ have isolated singularities only, and  set $p=mult(V,0)$.  Then:
\beq\label{Eq Milnor Number For Blowup}
\mu(V,0)=(p-1)^n+\mu(\mP^{n-1},\mP T_{(V,0)})+\mu(Bl_0\mC^n,\tV).
\eeq
As a preparation consider the change of the \ND~of $V_\al$ under the blowup of $(\mC^n,0)$.
\\\parbox{11cm}
{As in \S\ref{Sec tNnd directional approximation}, the top-dimensional faces of the diagram are
the ``face far from $\hat{x}_n$", denoted $\De$, and a collection $\{T_j\}$.
The intersections
$T_j\cap \De=t_j$ are faces of dimension $\le(n-2)$.
Consider the strict transform of $V_\al$ under the blowup of the origin
$\{(x_1,\cdots,x_n)=[\si_1,\cdots,\si_n]\}\subset\mC^n\times\mP^{n-1}$.
(The two \ND s are on the right.)
The relevant chart is $\si_n\neq0$, with the coordinates $(\frac{\si_1}{\si_n},\ldots,\frac{\si_{n-1}}{\si_n},x_n)$.
The total transform of the function is:\\
$f(x_1,\cdots,x_n)\to x^p_n\Big(f_p(\frac{\si_1}{\si_n},\cdots,\frac{\si_{n-1}}{\si_n},1)+
x_nf_{p+1}(\frac{\si_1}{\si_n},\cdots,\frac{\si_{n-1}}{\si_n},1)+\cdots+x^{q-p}_n\Big)+\cdots
$\\
The set of points under $\Ga_{(V_\al,0)}$ is naturally subdivided into two parts, one being the pyramid under
$Span(x^p_1,\ldots,x^p_n)$. Denote the other part (under $\cupl_j T_j$ but not under $Span(x^p_1,\ldots,x^p_n)$)
by $\Box_{V_\al}$. Let $\tBox_{V_\al}$ be the polyhedron under $\Ga_{(\tV_\al,0)}$.
}
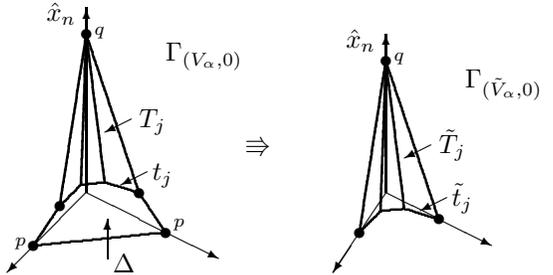
\begin{picture}(0,0)(-30,20)
\put(0,0){\vector(0,1){70}}\put(0,0){\vector(2,-1){50}}\put(0,0){\vector(-1,-1){30}}
\put(-22.5,-22.5){$\bullet$}\put(27.5,-17.5){$\bullet$}\put(-12.5,-7.5){$\bullet$}\put(17.5,-2.5){$\bullet$}

\put(-2.5,57.5){$\bullet$}  \put(3,60){$\tinyA q$}

\put(-15,65){$\hat{x}_n$}

\thicklines
\curve(-20,-20,-10,-5)\curve(30,-15,20,0)\curve(30,-15,-20,-20)
\curve(-10,-5,-2,3)\curve(7,4,-2,3)\curve(7,4,20,0)
\curve(0,60,-10,-5)\curve(0,60,20,0)\curve(0,60,7,4)\curve(0,60,-2,3)
\thinlines

\put(10,-30){$\De$}\put(8,-25){\vector(0,1){15}}
\put(25,5){$t_j$}\put(23,8){\vector(-2,-1){10}}  \put(20,25){$T_j$}\put(17,28){\vector(-2,-1){10}}
\put(-27,-20){$\tinyA p$}\put(33,-13){$\tinyA p$}  \put(30,50){$ \Ga_{(V_\al,0)}$}  \put(60,15){$\Rrightarrow$}
\end{picture}
\begin{picture}(0,0)(-140,30)
\put(0,10){\vector(0,1){60}}\put(0,10){\vector(2,-1){40}}\put(0,10){\vector(-2,-3){20}}
\put(-12.5,-7.5){$\bullet$}\put(17.5,-2.5){$\bullet$}

\put(-15,65){$\hat{x}_n$}

\put(-2.5,57.5){$\bullet$}  \put(3,60){$\tinyA q$}
\thicklines
\curve(-10,-5,-2,3)\curve(7,4,-2,3)\curve(7,4,20,0)
\curve(0,60,-10,-5)\curve(0,60,20,0)\curve(0,60,7,4)\curve(0,60,-2,3)
\thinlines

\put(25,5){$\tt_j$}\put(23,8){\vector(-2,-1){10}}  \put(20,25){$\tT_j$}\put(17,28){\vector(-2,-1){10}}
 \put(30,50){$ \Ga_{(\tV_\al,0)}$}
\end{picture}
\\\\
The blowup induces the map $\Box_{V_\al}\to\tBox_{V_\al}$
by $(a_1,\cdots,a_n)\to(a_1,\cdots,a_{n-1},\sum^n_{i=1}a_i-p)$. (Note that by construction $\sum a_i\ge p$.)
This map is an element of $SL(2,\mZ)$.
So,
there is a natural correspondence between the faces of $\Box_{V_\al}$ and $\tBox_{V_\al}$.
\\
While, $\Box_{V_\al}$ and $\tBox_{V_\al}$ are not equal as polyhedra, their corresponding faces have equal volumes.
Therefore $\mu(\Box_{V_\al})=\mu(\tBox_{V_\al})$. Here $\mu(\Box_{V_\al})$ is the standard expression: the main volume
$n!Vol_{n}(\Box_{V_\al})$ minus the volume of top-dimensional faces
$(n-1)!Vol_{n-1}(\Box_{V_\al})$, plus the sum of volumes of faces of codimension 2, etc.

Similarly, let $\Box^b_{V_\al}=\Box_{V_\al}\cap Span(x^p_1,\cdots,x^p_n)$
and  $\tBox^b_{V_\al}=\tBox_{V_\al}\cap\{x_n=0\}$.
Then the corresponding faces of $\Box^b_{V_\al}$ and $\tBox^b_{V_\al}$ have equal volume.

Therefore the classical Kouchnirenko formula gives
(recall that $\tV_\al$ and $\mP T_{(V_\al,0)}$ have \Nnd~singularities):
\beq
\mu(\Box_{V_\al})+\mu(\Box^b_{V_\al})=\mu(\tBox_{V_\al})+\mu(\tBox^b_{V_\al})=\mu(\tV_\al)+\mu(\mP T_{(V_\al,0)})
\eeq
Finally,
\li note that $\mu(\tV_\al)=\mu(\tV,pt_\al)$ and $\mu(\mP T_{(V_\al,0)})=\mu(\mP T_{(V,0)},pt_\al)$
\li sum over all the singular point of the projectivized tangent cone $\mP T_{(V,0)}$
(i.e. sum over all the directional approximations $V_\al$) and apply equation (\ref{Eq Milnor Number For Blowup})
\li note that each $\mu(V_\al,0)$ contains a contribution from the basic pyramid $(x^p_1,\cdots,x^p_n)$,
for which $\mu=(p-1)^n$.
\epr
\subsubsection{Proof of Proposition \ref{Thm Zeta Function for tNnd}}
Blowup $\mC^n$ at the origin. By the assumption the strict transform $\tV$ of $V$ has isolated
singularities only (e.g. $(\tV,y_1),\ldots,(\tV,y_k)$) and the exceptional divisor is $pE$
for $E\approx\mP^{n-1}\subset Bl_0(\mC^n)$. Now,
resolve the singularities of $\tV$.
\\
\parbox{14cm}
{Write the total preimage of the origin in the form:
$\pi^{-1}(0)=p\tE+\sum\limits_{\al=1}^kD_\al$.
Here each $D_\al$ corresponds to the resolution of $(\tV,y_\al)$ (cf. the picture).
In particular: $D_\al\cap D_\be=\empty$ for $\al\neq \be$.

Thus the product $\prod_{m\ge1}$ in the original formula (\ref{Eq Monodromy ACampo}) can
be replaced by $k$ copies (for each directional approximation $V_\al$).
Each such copy contributes the unnecessary factor \mbox{$(1-z^p)^{\chi(\tE\setminus(D_\al\cup\tV_\al))}$}
and no copy contains the needed factor \mbox{$(1-z^p)^{\chi(\tE\setminus(\cupl_\al D_\al\cup\tV))}$}.
}
\begin{picture}(0,0)(-20,0)
\curve(0,40,0,-30) \put(3,-30){$p\tE$}  \put(-2,28){$\bullet$}  \put(-2,0){$\bullet$}

\curve(0,30,20,35)\curve(10,35,30,25)\curve(20,25,40,35)
\curve(0,0,20,5)\curve(10,5,30,-5)\curve(20,-5,40,5)

\curvedashes{3,2}
\curve(0,35,30,40,50,35,55,35,60,30,55,25,50,25,30,20,0,25) \put(62,28){$D_\al$}
\curve(0,5,30,10,50,5,55,5,60,0,55,-5,50,-5,30,-10,0,-5) \put(62,-8){$D_\be$}
\end{picture}
\\\\
Hence, the formula can be written in the form
\beq
\ze_{(V,0)}(z)=(1-z^p)^{\chi(\tE\setminus(\cupl_\al D_\al\cup\tV))}.
\prod\limits_{\al=1}^k\frac{\ze_{(V_\al,0)}(z)}{(1-z^p)^{\chi(\tE\setminus(D_\al\cup\tV_\al))}}
\eeq
Note that (set-theoretically) $\tE\setminus(\cupl_\al D_\al\cup\tV)=E\setminus\mP T_{(V,0)}$ and
$\tE\setminus(D_\al\cup\tV_\al)=E\setminus\mP T_{(V_\al,0)}$, so the correction factor is:
\beq
\frac{(1-z^p)^{\chi(E\setminus\mP T_V)}}{\prod_\al(1-z^p)^{\chi(E\setminus\mP T_{V_\al})}}=
\frac{(1-z^p)^{(1-k)\chi(\mP^{n-1})}}{(1-z^p)^{\chi(\mP T_V)-\sum\chi(\mP T_{V_\al})}}
\eeq
Finally, as $\mP T_{(V_\al,0)}$ has isolated singularities only, $\chi(\mP T_{V_\al})=\chi(V_{p,n-1})+N_\al$,
where $V_{p,n-1}\subset\mP^{n-1}$ is a smooth hypersurface of degree $p$
(in particular its Euler characteristic is independent of the hypersurface) and $N_\al$ is
a number completely determined by the local singularity types of $Sing(\mP T_{(V_\al,0)})$.
Thus $\chi(\mP T_V)-\sum\chi(\mP T_{V_\al})=(1-k)\chi(V_{p,n-1})$, proving the statement.
\epr

\subsection{Some open questions}
\li The right modality for hypersurfaces. Already for surfaces there exist equisingular deformations involving monomials
below the \ND~\cite{Altmann91}. Hence a natural generalization of \S\ref{Sec.Right Modality in two variables} can possibly
give only a lower bound.

A candidate can be guessed from the case of Brieskorn-Pham singularities as follows. For the hypersurface $\{\sum x^{a_i}_i=0\}\subset(\mC^n,0)$
the miniversal equisingular deformation is spanned by monomials of the type:
\beq
\{x^{k_1}_1\cdots x^{k_n}_n|~~k_i\le a_i-2,~~\sum \frac{k_i}{a_i}\ge1\}.
\eeq
Apply to $\mZ^n$ the reflection: $(k_1,\cdots,k_n)\to(a_1-k_1,\cdots,a_n-k_n)$.
Then the relevant monomials are mapped to the region:
\beq
\{x^{k_1}_1\cdots x^{k_n}_n|~~2\le k_i\le a_i,~~\text{ the point }
(\frac{k_1}{n-1},\cdots,\frac{k_n}{n-1})\text{ is not above }\Ga\}.
\eeq
Note that for $n=2$ we obtain precisely the region of Kouchnirenko's formula
\ref{Sec Singularity Invariants The Case of Curves}.
Thus the natural guess is that for any $n$ and for any \Nnd~hypersurface whose \ND~intersects the
axes at the points $a_1,\cdots,a_n$,
the right modality is at least the number of these monomials.
One can easily check that this holds for suspensions of curves and $T_{pqr}$ singularities.

\li Recall that by \cite{Saito88},
for surfaces (i.e. $n=3$) the geometric genus of a Newton-non-degenerate isolated
singularity equals the number of $\mZ^n_{>0}$ not above $\Ga$.
Can this formula be generalized to the \dNnd~case?

\li For a \tNnd~singularity $f^{-1}(0)\sset(\mC^n,0)$ it is interesting to study the change of
the {\it analytical} type in the stratum
\beq
\Si_{f}:=\{g\in\mC\{x\} ~ | ~ \text{$\Ga_f=\Ga_g$ in any coordinate system}\}.
\eeq
While possible such $g$ are very restricted, e.g. the tangent cone is constant along this stratum,
some analytic invariants (e.g. the Tjurina number) can vary, see e.g.
the normal family of $E_{12}$: $x^3+y^7+axy^5$. A natural question is:
which analytic invariants are preserved in such families?

\li Can one distinguish the spectra of \Nnd, \dNnd, \tNnd~ from the spectra of other singularities?

\li Give an algorithm to check that a given hypersurface is \tNnd~ (or \dNnd~, or \Nnd~in some coordinates).


\begin{thebibliography}{99}
\bibitem[Abderrahmane05]{Abderrahmane05} O.Abderrahmane, {\it On the Lojasiewicz exponent and Newton polyhedron.}
 Kodai Math. J. 28 (2005), no. 1, 106--110
\bibitem[A'Campo75]{A'Campo75}N.A'Campo, {\it La fonction z\^{e}êta d'une monodromie.}
  Comment. Math. Helv.  50  (1975), 233--248

\bibitem[AGLV-book]{AGLV} V.I.Arnol'd, V.V.Goryunov, O.V.Lyashko, V.A.Vasil'ev, {\it Singularity theory. I.}
Reprint of the original English edition from the series Encyclopaedia of Mathematical Sciences
[ Dynamical systems. VI, Encyclopaedia Math. Sci., 6, Springer, Berlin, 1993]. Springer-Verlag, Berlin, 1998. iv+245 pp. ISBN: 3-540-63711-7

\bibitem[AGV-book]{AGV}V.I.Arnol'd, S.M.Gusein-Zade,A.N.Varchenko, {\it Singularities of differentiable maps. Vol. II.
Monodromy and asymptotics of integrals}. Monographs in Mathematics, 83. Birkhäuser Boston, Inc., Boston, MA, 1988.


\bibitem[Arnol'd74]{Arnol'd74}V.I. Arnol'd, {\it Normal forms of functions in the neighborhood of degenerate critical points.}
Uspehi Mat. Nauk  29  (1974), no. 2(176), 11--49

\bibitem[AB-L-MH06]{AB-L-MH06}E.Artal-Bartolo, I.Luengo, A.Melle-Hern\'{a}ández,  {\it Superisolated surface
singularities.} Singularities and computer algebra, 13--39, London Math. Soc. Lecture Note Ser., 324, Cambridge Univ.
Press, Cambridge, 2006

\bibitem[Altmann91]{Altmann91}K.Altmann, {\it Equisingular deformations below the Newton boundary.}
  Compositio Math.  80  (1991),  no. 3, 257--283


\bibitem[BGMM89]{BGMM89}J.Brian\c{c}çon, M.Granger, Ph.Maisonobe, M.Miniconi,
{\it Algorithme de calcul du polyn\^{o}ôme de Bernstein: cas non d\'{e}ég\'eén\'eér\'eé.}Ann. Inst. Fourier (Grenoble)  39  (1989),  no. 3, 553--610.


\bibitem[Brian\c{c}on-Speder75]{Briancon-Speder75} J.Brian\c{c}çon, J.P.Speder, {\it La trivialit\'{e}é
topologique n'implique pas les conditions de Whitney. }  C. R. Acad. Sci. Paris S\'{e}ér. A-B  280
(1975), no. 6, Aiii, A365--A367

\bibitem[Danilov-Khovanski87]{Danilov-Khovanski87}V.I.Danilov, A.G.Khovanski\u{i},
{\it Newton polyhedra and an algorithm for calculating Hodge-Deligne numbers.}
Izv. Akad. Nauk SSSR Ser. Mat. 50 (1986), no. 5, 925--945

\bibitem[Fukui91]{Fukui91} T.Fukui, {\L}ojasiewicz type inequalities and Newton diagrams. Proc.
Amer. Math. Soc. 112 (1991), no. 4, 1169--1183

\bibitem[GaLeP{\l}07]{GaLePl07}E.Garc\'{\i}a Barroso, A.Lenarcik, A.P{\l}oski{\it Characterization of non-degenerate plane curve singularities}
preprint, arXiv:0711.2833


\bibitem[Gusein-Zade Luengo Melle-Hernandez97]{Gusein-Zade Luengo Melle-Hernandez97} S.M.Gusein-Zade,
I.Luengo, A.Melle-Hern\'{a}ández,
{\it Partial resolutions and the zeta-function of a singularity}.  Comment. Math. Helv.  72  (1997),  no. 2, 244--256



\bibitem[GLS-book]{GLSbook} G.-M. Greuel, C. Lossen, E. Shustin:
{\it Introduction to Singularities and Deformations}.
Series: Springer Monographs in Mathematics 2006. ISBN: 3-540-28380-3

\bibitem[GPS-Singular]{Singular} G.M.Greuel, G.Pfister, H.Sch\"onemann, {\it {\sc Singular} 3.0.0.
A Computer Algebra System for Polynomial Computations.} ZCA, University of
Kaiserslautern (2005). {http://www.singular.uni-kl.de}.



\bibitem[Kerner06]{Kerner06}D.Kerner, {\it Enumeration of singular algebraic curves.}
Israel J. Math. 155 (2006), 1--56.


\bibitem[Kerner-OWP]{Kerner-OWP}D.Kerner (OWLF 2008)
{\it On the Directionally Newton-non-degenerate Singularities of Complex Hypersurfaces}
Oberwolfach preprint 2008 - 16, http://www.mfo.de/publications/owp/


\bibitem[Khovanski-Varchenko85]{Khovanski-Varchenko85}A.G.Khovanski\u{i},A.N.Varchenko,
{\it Asymptotic behavior of integrals over vanishing cycles and the Newton polyhedron.}
Dokl. Akad. Nauk SSSR  283 (1985),  no. 3, 521--525

\bibitem[Kouchnirenko76]{Kouchnirenko76}A.G.Kouchnirenko, {\it Poly\'{e}èdres de Newton et nombres de Milnor.}
Invent. Math. 32 (1976), no. 1, 1--31


\bibitem[Kulikov98]{Kulikov98}V.S.Kulikov, {\it Mixed Hodge structures and singularities.} Cambridge Tracts
in Mathematics, 132. Cambridge University Press, Cambridge, 1998. xxii+186 pp

\bibitem[L{\^{e}}-Ramanujam76]{Le-Ramanujam76}D.T.L{\^{e}} C.P. Ramanujam, {\it The invariance of Milnor's number
implies the invariance of the topological type.} Amer. J. Math. 98 (1976), no. 1, 67--78

\bibitem[Luengo87]{Luengo87}I.Luengo, {\it The $µ$-constant stratum is not smooth.}
Invent. Math. 90 (1987), no. 1, 139--152


\bibitem[Luengo-Melle95]{Luengo-Melle95}I.Luengo, A.Melle, {\it A formula for the Milnor number.}
C. R. Acad. Sci. Paris S\'{e}ér. I Math. 321 (1995), no. 11, 1473--1478


\bibitem[Melle00]{Melle00}A.Melle-Hernández, {\it Milnor numbers for surface singularities.}
  Israel J. Math.  115  (2000), 29--50.


\bibitem[Oka79]{Oka79}M.Oka, {\it On the bifurcation of the multiplicity and topology of the Newton boundary}.
J. Math. Soc. Japan 31 (1979), no. 3, 435--450

\bibitem[Parusin'ski99]{Parusin'ski99}A.Parusin'ski, {\it Topological triviality of $µ$-constant deformations of
type $f(x)+tg(x)$.} Bull. London Math. Soc. 31 (1999), no. 6, 686--692


\bibitem[Saito88]{Saito88}M.Saito, {\it Exponents and Newton polyhedra of isolated hypersurface singularities.}
Math. Ann.  281  (1988),  no. 3, 411--417.


\bibitem[Steenbrink76]{Steenbrink76}J.M.Steenbrink, {\it Mixed Hodge structure on the vanishing cohomology}.
Real and complex singularities (Proc. Ninth Nordic Summer School/NAVF Sympos. Math., Oslo, 1976),
 pp. 525--563. Sijthoff and Noordhoff, Alphen aan den Rijn, 1977


\bibitem[Steenbrink85]{Steenbrink85}J.H.M.Steenbrink, {\it Semicontinuity of the singularity spectrum.}
  Invent. Math.  79  (1985),  no. 3, 557--565.

\bibitem[Teissier77]{Teissier77}B.Teissier, {\it Vari\'{e}t\'{e}s polaires. I.
Invariants polaires des singularit\'{e}s d'hypersurfaces.}
Invent. Math.  40  (1977), no. 3, 267--292.


\bibitem[Varchenko76]{Varchenko76}A.N.Varchenko, {\it Zeta-function of monodromy and Newton's diagram.}
Invent. Math. 37 (1976), no. 3, 253--262

\bibitem[Varchenko82]{Varchenko82}A.N.Varchenko, {\it A lower bound for the codimension of the $\mu={\rm const}$
stratum in terms of the mixed Hodge structure.} Vestnik Moskov. Univ. Ser. I Mat. Mekh.  1982,  no. 6, 28--31, 120


\end{thebibliography}
\end{document}